\documentclass[a4paper]{article}

\usepackage{a4wide}
\usepackage{latexsym}
\usepackage[UKenglish]{babel}
\usepackage{graphicx,subfigure}
\usepackage{algorithm}
\usepackage{algpseudocode}
\usepackage{epstopdf}
	\graphicspath{{./}{figures/}}
\usepackage[colorlinks=true,linkcolor=blue,citecolor=blue]{hyperref}
\usepackage{xcolor}
\usepackage{amsfonts,amsmath,amsthm,mathtools,bbm,cool}

\newtheorem{theorem}{Theorem}[section]
\newtheorem{proposition}[theorem]{Proposition}
\theoremstyle{remark}\newtheorem{remark}[theorem]{Remark}

\newcommand{\be}{\begin{equation}}
\newcommand{\ee}{\end{equation}}

\allowdisplaybreaks

\begin{document}
\title{Kinetic-controlled hydrodynamics for multilane traffic models}

\author{Raul Borsche \\
		{\small Department of Mathematics	} \\
		{\small TU Kaiserslautern, Germany} \\
		{\small\tt borsche@mathematik.uni-kl.de } \\[5mm]	
		Axel Klar \\
		{\small Department of Mathematics	} \\
		{\small TU Kaiserslautern, Germany} \\
		{\small\tt klar@mathematik.uni-kl.de } \\[5mm]
		Mattia Zanella \\
		{\small	Department Mathematics ``F. Casorati''} \\
		{\small University of Pavia, Italy} \\
		{\small\tt mattia.zanella@unipv.it}}
\date{}

\maketitle

\begin{abstract}
We study the application of a recently introduced hierarchical description of traffic flow control by driver-assist vehicles to include lane changing dynamics. Lane-dependent feedback control strategies are implemented at the level of vehicles and the aggregate trends are studied by means of Boltzmann-type equations determining three different hydrodynamics based on the lane switching frequency. System of first order macroscopic equations describing the evolution of densities along the lanes are then consistently determined through a suitable closure strategy. Numerical examples are then presented to illustrate the features of the proposed hierarchical approach.
\medskip


\end{abstract}


\section{Introduction}
\label{sect:intro}

In recent years an increased level of automation invested transportation industry thanks to new technologies in algorithmic sensing. Even if the automation of vehicular traffic is nowadays at its germinal state, the rapid changes that we are facing in modern societies in the definition of sustainable management of resources should inspire substantial questions about traffic management, its governance and its relations with a zero-emission society. All these aspects can boost effective  developments in driver-assist and self-driving cars, having a tremendous impact in designing new transportation systems \cite{Chow}. Furthermore, among the main goals of such technologies is the enhancement of driver's safety through the mitigation of road risk factors, like the ones related to the heterogeneity of velocities of the traffic stream, that are responsible, as reported in \cite{WHO}, for an increase of crash risk, see also the field experiment \cite{Stern}.

Recent efforts have been devoted to understand the impact of automation on observable macroscopic quantities characterizing the vehicular flows, see e.g. \cite{DNP,GGLP}. Among them we highlight the role of kinetic modelling of vehicular traffic which is capable to organically link microscopic dynamics, where autonomous decision-making processes are active, with macroscopic emergent features of traffic \cite{ChPiTo,PTZ}. The main idea relies in incorporating optimal control strategies in the definition of speeds' variations \cite{TZ0,TZ1} whose observable effects are then consistently defined. 

In the present work we formalize a kinetic multi-lane traffic flow model where driver-assist controls are activated in a portion of the traffic stream allowing vehicles to switch lane if the local density in the neighboring lane is sufficiently small. The idea in this context is to exploit the possibility to control few vehicles, i.e. the automated cars, to reduce speed dependent risk factors and to induce a regularization of the traffic flow. Recent experimental observation highlighted the emergence of Beta-type speed distributions \cite{HTVZ2,NHJ}. For this reason, the in-lane dynamics is supposed to be given by a kinetic model based on anisotropic binary interactions recently introduced in \cite{HTVZ2,TZ1,TZ2} and reproducing the class of observed speed distributions and that are consistent with macroscopic features related to the fundamental diagram of traffic. Other binary interaction models producing consistent fundamental diagrams of traffic have been designed in the literature \cite{HP,HTVZ,VHPT}.  
 Regarding the dynamics between lanes, without intending to review the huge literature on this topic, we highlight that the designed switching dynamics has been introduced in agreement with existing modeling for traffic dynamics at the mean-field \cite{GPV}, kinetic scale \cite{HTVZ,IKM,KW1,KW2} and at the macroscopic scale \cite{HFV,Song}. The importance of incorporating lane switching is fundamental in traffic modelling since it incorporates an important realistic feature of the traffic phenomena under study which should be taken into account for a robust implementation of driver-assist technologies. 

The description and modelling of emerging phenomena in traffic dynamics is related to the possibility to describe the effects at different scales of possible external control actions \cite{PTZ}. To this end we will define suitable controls to modify speed-based costs whose impact essentially depend on the switching dynamics. Indeed we will show how, in relation with switching coefficients, radically different hydrodynamics, i.e. macroscopic models describing evolution of the lane densities, may be computed starting from the kinetic-controlled model. In details three main regimes are detected: the collision dominated regime, when in-lane interaction dynamics are prevalent with respect to lane crossing, the fast switching regime, when viceversa lane crossing is more frequent than drivers' interactions, and finally the slow switching regime, emerging as a precise balance between the two possible behaviors. It will be proven how the introduced control is capable to modify the flux of first order hydrodynamic equations for traffic flow by inducing alignment of in-lane speeds. 

In more details, the paper is structured as follows. In Section \ref{sect:multilane} we design a multilane kinetic model where binary control problems are introduced to mimic the action of driver-assist vehicles taking into account their penetration rate into the traffic stream. Hence, we introduce an inhomogeneous Boltzmann-type model which incorporates also lane changing which takes place locally in space and is based on both the densities of the two lanes. An analytic insight is given to determine the evolution of densities, mean speeds and energies of the lanes. In particular, the fundamental diagram of the controlled dynamics is studied. In Section \ref{sec:hydro} we derive three different hydrodynamic models for the constrained traffic dynamics by distinguishing between three main regimes for  lane changing. Finally, in Section \ref{sec:numerics} we investigate the solution of the kinetic and hydrodynamic models by means of suitable numerical methods to observe the impact of the introduced driver-assist technologies.  

\section{Multilane kinetic model with driver-assist controls }\label{sect:multilane}
Let us consider a two lanes road. We introduce therefore two distribution functions $f_1 = f_1(x,v,t)$ and $f_2 = f_2(x,v,t)$ such that $f_1\,dv\,dx$ and $f_2\,dv\,dx$ are the fraction of vehicles at time $t>0$  in the first and second lane, respectively, and traveling with speeds in $[v,v+dv]$ in $[x,x+dx]$. All the variables are dimensionless and, in particular, we consider normalized speeds $v\in[0,1]$. 

As for classical kinetic theory, we define the densities $\rho_1 = \rho_1(x,t)$ and $\rho_2 = \rho_2(x,t)$ as follows
\[
\rho_1(x,t) = \int_{0}^1 f_1(x,v,t)dv,\qquad \rho_2(x,t) = \int_{0}^1 f_2(x,v,t)dv.
\]
In a Boltzmann-type kinetic approach for traffic dynamics we study the evolution of the introduced distributions by means of microscopic stochastic binary interactions among drivers, which determine the registered speed change. Suitable driver-assist controls will be implemented at the vehicles' level and the upscaled. 

\subsection{Constrained binary interactions}
If we denote the pre-interaction speeds for the vehicles in the first lane $(v_1,w_1)$ and $(v_2,w_2)$ for the vehicles in the second lane, we may express their post-interaction states $(v_1^\prime,w_1^\prime)$ and $(v_2^\prime,w_2^\prime)$ using the formalism
\begin{equation}
\label{eq:v1-2}
\begin{split}
&v_i^\prime = v_i + \gamma\left( I(v_i,w_i,\rho_i) + \Theta \mathfrak{u}_i \right) + D(v_i,\rho_i)\eta_i,\\
&w_i^\prime = w_i,
\end{split}
\end{equation}
with $i =1,2$.
In \eqref{eq:v1-2} we denoted with $\gamma>0$ a proportionality parameter and with $I(\cdot,\cdot,\cdot)$ a general interaction rule expressing the speed change through interaction with another vehicle in a given density regime. Centered random variables $\eta_1$, $\eta_2$ with finite variances $\sigma_1^2$ and $\sigma_2^2$ have been introduced, expressing  the deviations from the prescribed dynamics, while the function $D(\cdot,\cdot)>0$ is a  local  diffusion coefficient. \\

It is easily observed how a consistent part of the modeling efforts for traffic dynamics relies in the definition of the interaction rule. Without the aim to be exhaustive we point the interested reader to \cite{HP,HTVZ,VHPT} and the references therein for several examples. More recently in \cite{TZ1} the following interaction rule has been introduced
\begin{equation}
\label{eq:inter}
I(v_i,w_i,\rho_i) = P(\rho_i) (1-v_i) + (1-P(\rho_i))(P(\rho_i)w_i-v_i), \qquad i = 1,2. 
\end{equation}

This rule expresses the tendency of an agent to naturally relax towards the maximal nondimensional allowed speed and, at the same time, to decelerate to a portion of the speed of the vehicle with which the interaction takes place. This behavior is expressed in \eqref{eq:inter} as a convex combination weighted by the dimensionless density-dependent function $P(\cdot)$ whose general form can be assumed as 
\begin{equation}\label{eq:defP}
P(\rho_i) = (1-\rho_i)^{\mu}, \qquad \mu >0. 
\end{equation}
The function $P(\cdot)$,  expressing the probability of acceleration, is typically linked to the class of vehicles of the traffic flow, see \cite{HTVZ2}. In the present work we fix for simplicity a deterministic $\mu>0$ or, in other words, we concentrate on a traffic flow composed by a single class of vehicles. 

The terms $\mathfrak{u}_i$ are controls representing instantaneous corrections of the usual interactions operated by  vehicles equipped with a driver assist technology. Since in real traffic dynamics the number of those automated cars is usually small, we multiplied the controls by a Bernoulli random variable $\Theta\in\{0,1\}$ with
that $\Theta \sim \textrm{Bernoulli}(p)$, where  $p \in[0,1]$ is a given probability,which 
can be seen as the fraction of driver-assist vehicles in the traffic stream, i.e. the penetration rate. 

The controls $\mathfrak{u}_i$, $i=1,2$ are chosen in such a way that  they optimize a suitable binary cost functional $J = J(v^\prime,\mathfrak{u})$, whose minimization links the control to some recommended behavior, i.e.
\begin{equation}
\label{eq:argmin}
\mathfrak{u}_i = \textrm{arg}\min_{\mathfrak{u}_i\in \mathcal U} J_i(v_i^\prime, \mathfrak u_i),\qquad i = 1,2,
\end{equation}
subject to \eqref{eq:v1-2},where   $\mathcal U$ is the set of admissible controls such that $v_i^\prime \in [0,1]$. Among the possible cost functionals we will consider in the following 
\[
J_i (v_i,\mathfrak{u}_i) = \dfrac{1}{2} \mathbb E_{\eta_i} \left[ (\bar v_i(\rho_i)-v_i^\prime)^2 + \nu_i \mathfrak{u}^2 \right],\qquad i = 1,2,
\]
where $\mathbb E_{\eta_i}[\cdot]$ denotes the expectation with respect to the random variable $\eta_i$, so that the control mechanism minimize the $L^2$ distance of the post-interaction speeds with some recommended speed $\bar v_i \in[0,1]$ which may depend on the lane's density, i.e. $\bar v_i = \bar v_i(\rho_i)$. 

\subsection{Optimal binary controls}

We may determine the minimization \eqref{eq:argmin} subject to  the dynamics \eqref{eq:v1-2}  through a standard Lagrange multiplier approach. Let us define the Lagrangian
\[
\mathcal L(v^\prime_i,\mathfrak{u}_i,\lambda) = J_i(v^\prime_i,\mathfrak{u}_i) + \lambda \mathbb E_{\eta_i}\left[ v_i^\prime - v_i - \gamma \left( I(v_i,w_i,\rho_i) + \Theta \mathfrak u_i \right) - D(v_i,\rho_i)\eta \right],
\]
being $\lambda \in \mathbb R$ the Lagrange multiplier associated to the constraint \eqref{eq:argmin}. Hence, we can compute
\begin{equation}
\begin{cases}
\partial_{\mathfrak u} \mathcal L = \nu_i \mathfrak u_i - \gamma \Theta \lambda \\
\partial_{v_i^\prime} \mathcal L = \left \langle v_i^\prime - \bar v_i(\rho_i) \right\rangle + \lambda.
\end{cases}
\end{equation}
By imposing $\partial_{\mathfrak u} \mathcal L =0$ and $\partial_{v_i^\prime} \mathcal L =0$ we have the lane-dependent optimal control
\[
\mathfrak u^*_i = \dfrac{\gamma}{\nu_i } \mathbb E_{\eta_i}\left[ \bar v_i(\rho_i) - v_i^\prime \right], 
\]
which can be expressed as a function of the pre-interaction speeds  $(v_i,w_i)$ as follows
\begin{equation}
\mathfrak u^*_i = \dfrac{\gamma \Theta}{\nu_i + \gamma^2 \Theta^2} (\bar v_i(\rho_i) - v_i) - \dfrac{\gamma^2\Theta}{\nu_i + \gamma^2 \Theta^2} I(v_i,w_i,\rho_i). 
\end{equation}

We obtained a feedback formulation of the introduced control which can be inserted in the original binary interaction rules which modify then as follows
\begin{equation}
\label{eq:v1-2_control}
\begin{split}
&v_i^\prime = v_i + \dfrac{\nu_i \gamma}{\nu_i + \gamma^2 \Theta^2} I(v_i,w_i,\rho_i) + \dfrac{\gamma^2 \Theta^2}{\nu_i + \gamma^2 \Theta^2} (\bar v_i - v_i) + D(v_i,\rho_1)\eta_1,\\
&w_i^\prime = w_i,
\end{split}
\end{equation}
with $i=1,2$.
 With the following result  we can guarantee the physical admissibility of the control actions, see Proposition 3.1 in \cite{TZ1}. 
 \begin{proposition}
 In \eqref{eq:v1-2_control} let $\gamma\in[0,1]$, $\nu_i>0$ and $I(v_i,w_i,\rho_i)$ like in \eqref{eq:inter} for $i=1,2$. If there exists $c_i>0$ such that
 \[
 \begin{cases}
 |\eta_i| \le c_i\left(1-\dfrac{\nu_i + \gamma}{\nu_i}\gamma \right), \\
 c_i D(v_i,\rho_i) \le  \min\{v_i,1-v_i\}
 \end{cases}
 \]
then the pre-interaction speeds $v_i^\prime \in [0,1]$ for any $v_i,w_i \in[0,1]$, $\rho_i \in [0,1]$ and with $i=1,2$. 
 \end{proposition}

\subsection{Inhomogeneous Boltzmann-type model for multilane dynamics }

Let us introduce a inhomogeneous Boltzmann-type equation modeling multilane traffic dynamics
\begin{equation}
\label{eq:multi}
\begin{split}
\partial_t f_1 + v\partial_x f_1 = Q(f_1,f_1) - \beta_1 (1-\rho_2)^\alpha f_1  +\beta_2 (1-\rho_1)^\alpha f_2, \\
\partial_t f_2 + v\partial_x f_2 = Q(f_2,f_2) - \beta_2 (1-\rho_1)^\alpha f_2  + \beta_1(1-\rho_2)^\alpha f_1,\\
\end{split}
\end{equation}
where  $\alpha>0$ measures the sensitivity of lane changing with respect to the density of the adjacent lane, and where $Q(f_i,f_i)$, $i=1,2$ is the Boltzmann collision operator with localized interactions
\begin{align}
\label{collision}
Q(f_i,f_i)(x,v,t) = \dfrac{1}{2} \mathbb E_{\Theta,\eta_i}\left[ \int_0^1 \left( \dfrac{1}{{}^\prime J_i} f_i(x,{}^\prime v,t)f_i(x,{}^\prime w,t) - f_i(x,v,t)f_i(x,w,t) \right)dw \right],
\end{align}
where $({}^\prime v,{}^\prime w)$ are the pre-interaction speed sampled from the distribution $f_i$ and generating the post-interaction speeds according to the binary rules \eqref{eq:v1-2_control} and ${}^\prime J_i$ is the Jacobian of the transformation. 
 
The additional reaction terms mimic a dynamics where a given vehicle belonging to lane $i$ changes its lane with probability $\beta_i (1-\rho_j)^\alpha$, $j \ne i$, where $\rho_j$ is the density of the second lane $j$. 
We point the interested reader to  \cite{IKM,KW1,KW2,Song} for a more detailed multi-lane model at the kinetic and macroscopic scale. 

In order to compute evolution of observable quantities, it is more convenient to recast \eqref{eq:multi} in weak form. Let us introduce a test function $\varphi(v)$, we have
\begin{equation}\label{eq:1_weak}
\begin{split}
\partial_t & \int_0^1 \varphi(v)f_1(x,v,t)dv +\partial_x \int_0^1 v\varphi(v) f_1(x,v,t)dv = \\
& \dfrac{1}{2}\mathbb E_{\Theta,\eta_i}\left[ \int_0^1 \int_0^1 (\varphi(v^\prime)-\varphi(v)) f_1(x,v,t)f_1(x,w,t)\,dv\,dw \right]  \\
& -\beta_1(1-\rho_2)^\alpha \int_0^1 \varphi(v)f_1(x,v,t)dv + \beta_2 (1-\rho_1)^\alpha \int_0^1 \varphi(v)f_2(x,v,t)dv,
\end{split}
\end{equation}
and 
\begin{equation}\label{eq:2_weak}
\begin{split}
\partial_t & \int_0^1 \varphi(v)f_2(x,v,t)dv +\partial_x \int_0^1 v\varphi(v) f_2(x,v,t)dv = \\
& \dfrac{1}{2}\mathbb E_{\Theta,\eta_i}\left[ \int_0^1 \int_0^1 (\varphi(v^\prime)-\varphi(v)) f_2(x,v,t)f_2(x,w,t)\,dv\,dw\right]  \\
& -\beta_2(1-\rho_1)^\alpha \int_0^1 \varphi(v)f_2(x,v,t)dv + \beta_1 (1-\rho_2)^\alpha \int_0^1 \varphi(v)f_1(x,v,t)dv.
\end{split}
\end{equation}

\subsubsection{Space homogeneous case: evolution of observable quantities}
In the space homogeneous setting, i.e. $f_1 = f_1(v,t)$ and $f_2 = f_2(v,t)$ we may compute the evolution of densities by setting $\varphi(v) = 1$ in \eqref{eq:1_weak}-\eqref{eq:2_weak}. We obtain the following nonlinear system of differential equations
\begin{equation}
\label{eq:system_rho}
\begin{cases}\vspace{0.2cm}
\dfrac{d}{dt} \rho_1 = -\beta_1 (1-\rho_2)^{\alpha} \rho_1 + \beta_2 (1-\rho_1)^{\alpha}\rho_2 = -\beta_1 \rho_1^{\alpha+1} + \beta_2 \rho_2^{\alpha+1} \\
\dfrac{d}{dt} \rho_2 = -\beta_2 (1-\rho_1)^{\alpha} \rho_2 + \beta_1 (1-\rho_2)^{\alpha}\rho_1 = -\beta_2 \rho_2^{\alpha+1} + \beta_1 \rho_1^{\alpha+1},
\end{cases}
\end{equation}
since in this case we may set $\rho_1(0) + \rho_2(0)= 1$, $\rho_i(0)\in [0,1]$, without loss of generality.  Furthermore, it is easily seen that the total mass $\rho_1+\rho_2$ is conserved in time for any choice of $\alpha>0$ and $\beta_1,\beta_2>0$. It is worth to notice that if $\beta_1 =\beta_2= 0$ no exchange is active and both densities are conserved. 
 
\begin{figure}
 \centering
 \includegraphics[scale = 0.35]{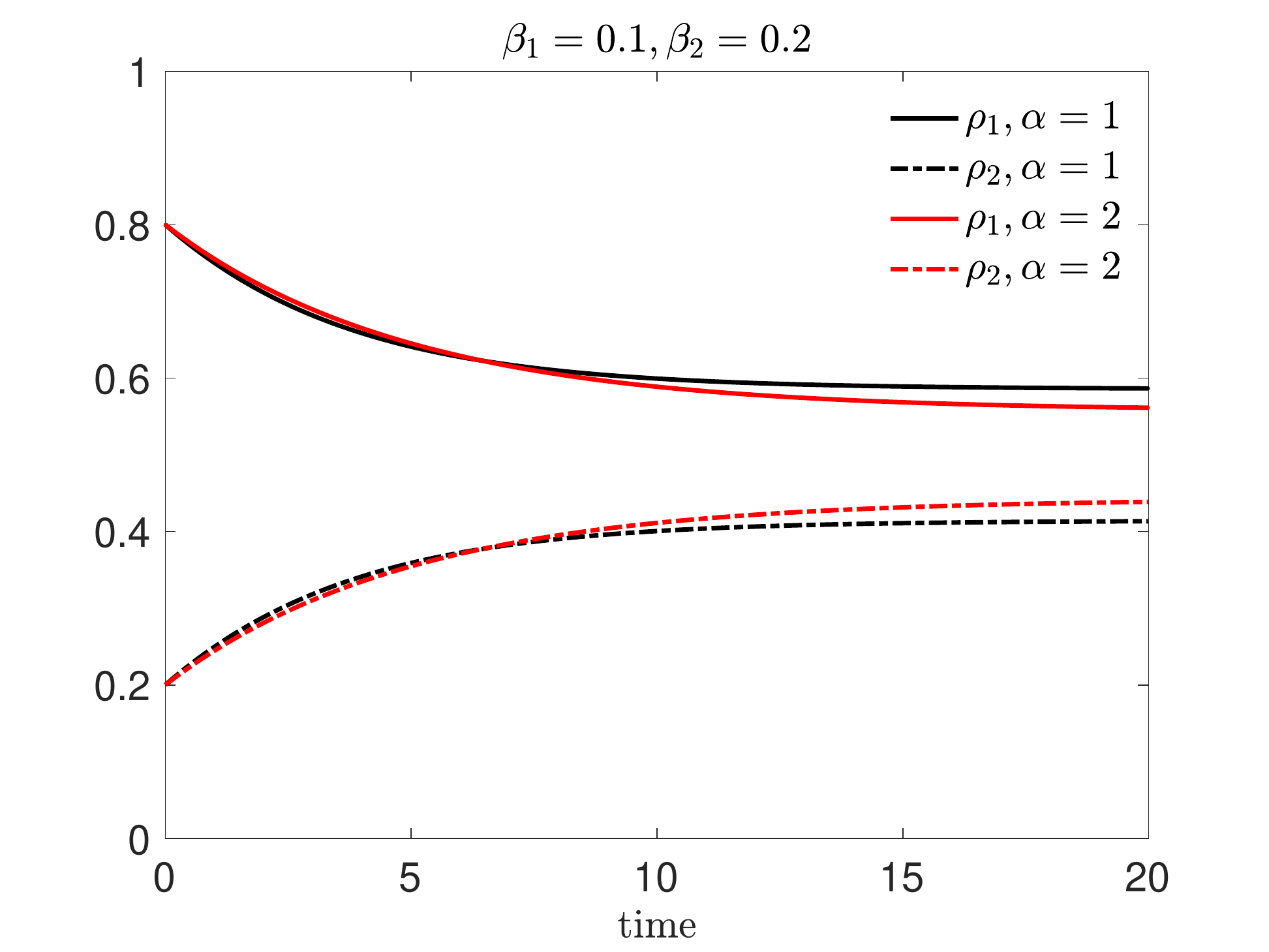}
  \includegraphics[scale = 0.35]{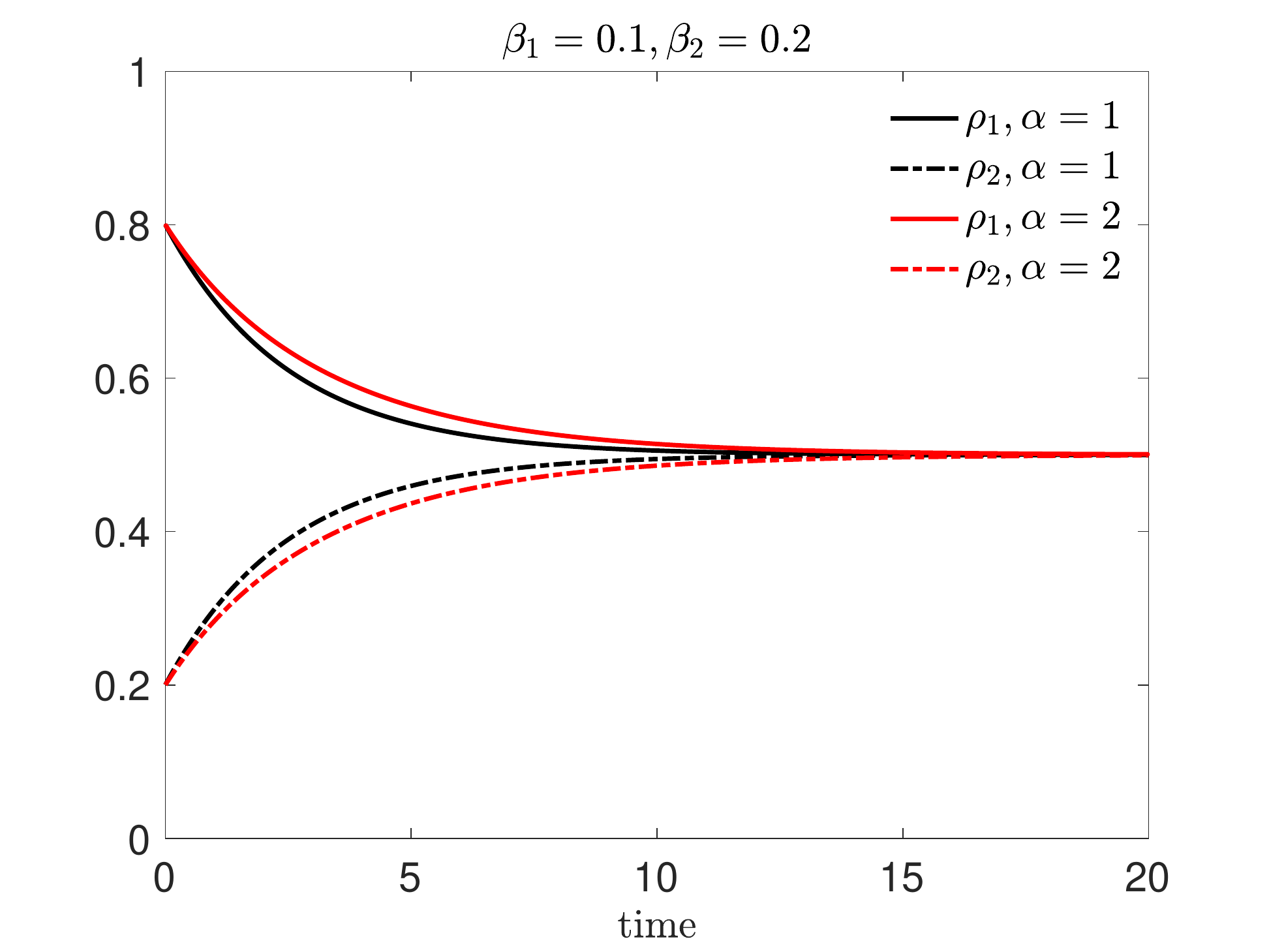}
  \caption{Evolution of the densities from \eqref{eq:system_rho} supposing $\rho_1 = 1-\rho_2 = 0.8$ and $\alpha =1$, $\alpha=2$. Left: asymmetric lane switching $\beta_1  = 0.1$ and $\beta_2 = 0.2$. Right: symmetric lane switching $\beta_1 = \beta_2 = 0.2$. }
  \label{fig:rho}
 \end{figure}

 In Figure \ref{fig:rho} we depict the evolution of the system \eqref{eq:system_rho} in the case $\rho_1(0) = 0.2$, $\rho_2(0) = 0.8$ and two possible values of the coefficients $\beta_1$, $\beta_2>0$. It is easily seen how the asymmetric case $\beta_1<\beta_2$ leads to a density reduction in the first lane and a density gain for the second lane. In the symmetric case $\beta_1= \beta_2=0.2$ we get an asymptotic balance of the two densities. 

Setting instead $\varphi(v) = v$ in \eqref{eq:1_weak}-\eqref{eq:2_weak} we get the evolution of the mean speed of the first and second lane, respectively
\[
m_1(t)= \dfrac{1}{\rho_1} \int_0^1 vf_1(v,t)dv,\qquad m_2(t) = \dfrac{1}{\rho_2 }\int_0^1 vf_2(v,t)dv.
\]
Their evolution is given by
\begin{equation}
\label{eq:m12_general}
\begin{cases}
\dfrac{d}{dt}\left(\rho_1 m_1 \right) =& \dfrac{\gamma \rho_1^2}{2} \left\{ \dfrac{\nu_1 + (1-p)\gamma^2}{\nu_1 + \gamma^2}\left( P(\rho_1) - (P(\rho_1)+(1-P(\rho_1))^2)m_1 \right) \right. \\
& \left. + \dfrac{\gamma p}{\nu_1 + \gamma^2 }  (\bar v_1 - m_1) \right\}  -\beta_1 (1-\rho_2)^{\alpha}\rho_1 m_1 + \beta_2 (1-\rho_1)^\alpha \rho_2m_2 \\ 
\dfrac{d}{dt}\left(\rho_2 m_2 \right) =& \dfrac{\gamma \rho_2^2}{2} \left\{ \dfrac{\nu_2 + (1-p)\gamma^2}{\nu_2 + \gamma^2}\left( P(\rho_2) - (P(\rho_2)+(1-P(\rho_2))^2)m_2 \right)  \right.\\ 
&\left.+ \dfrac{\gamma p}{\nu_2 + \gamma^2 }  (\bar v_2 - m_2) \right\}  -\beta_2 (1-\rho_1)^{\alpha}\rho_2 m_2 + \beta_1 (1-\rho_2)^\alpha \rho_1m_1.
\end{cases}
\end{equation}

\begin{figure}
\includegraphics[scale = 0.25]{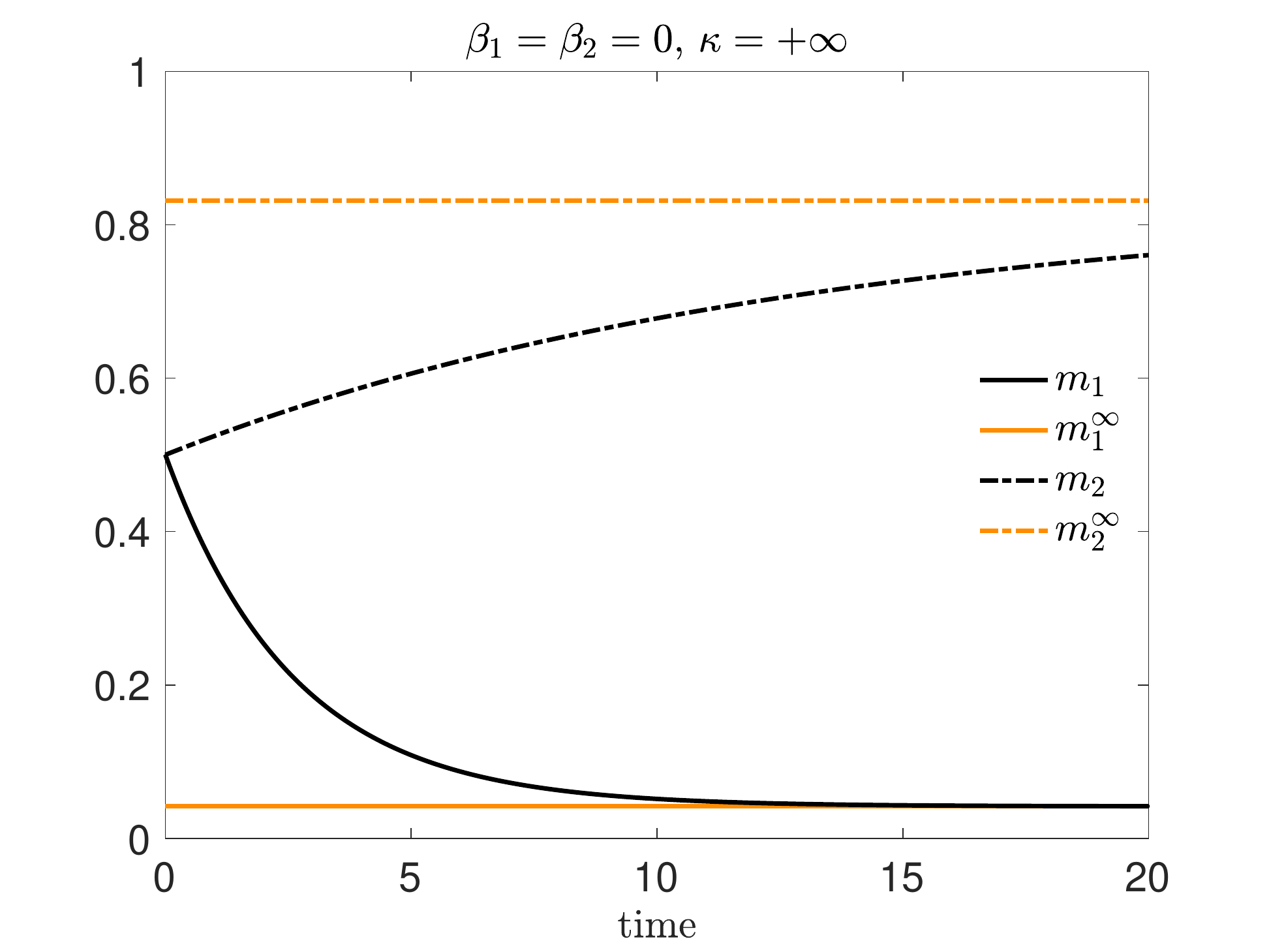}
\includegraphics[scale = 0.25]{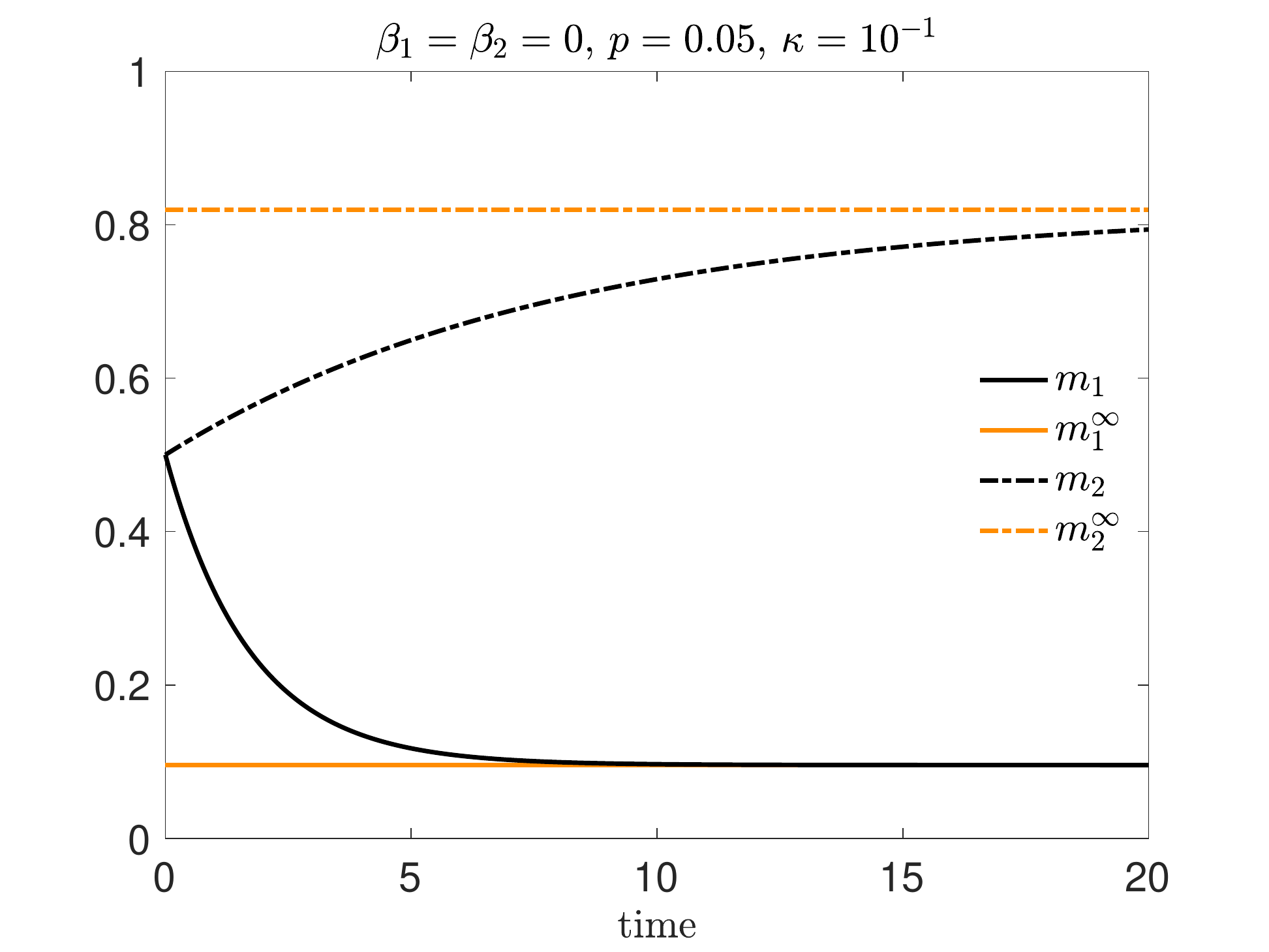}
\includegraphics[scale = 0.25]{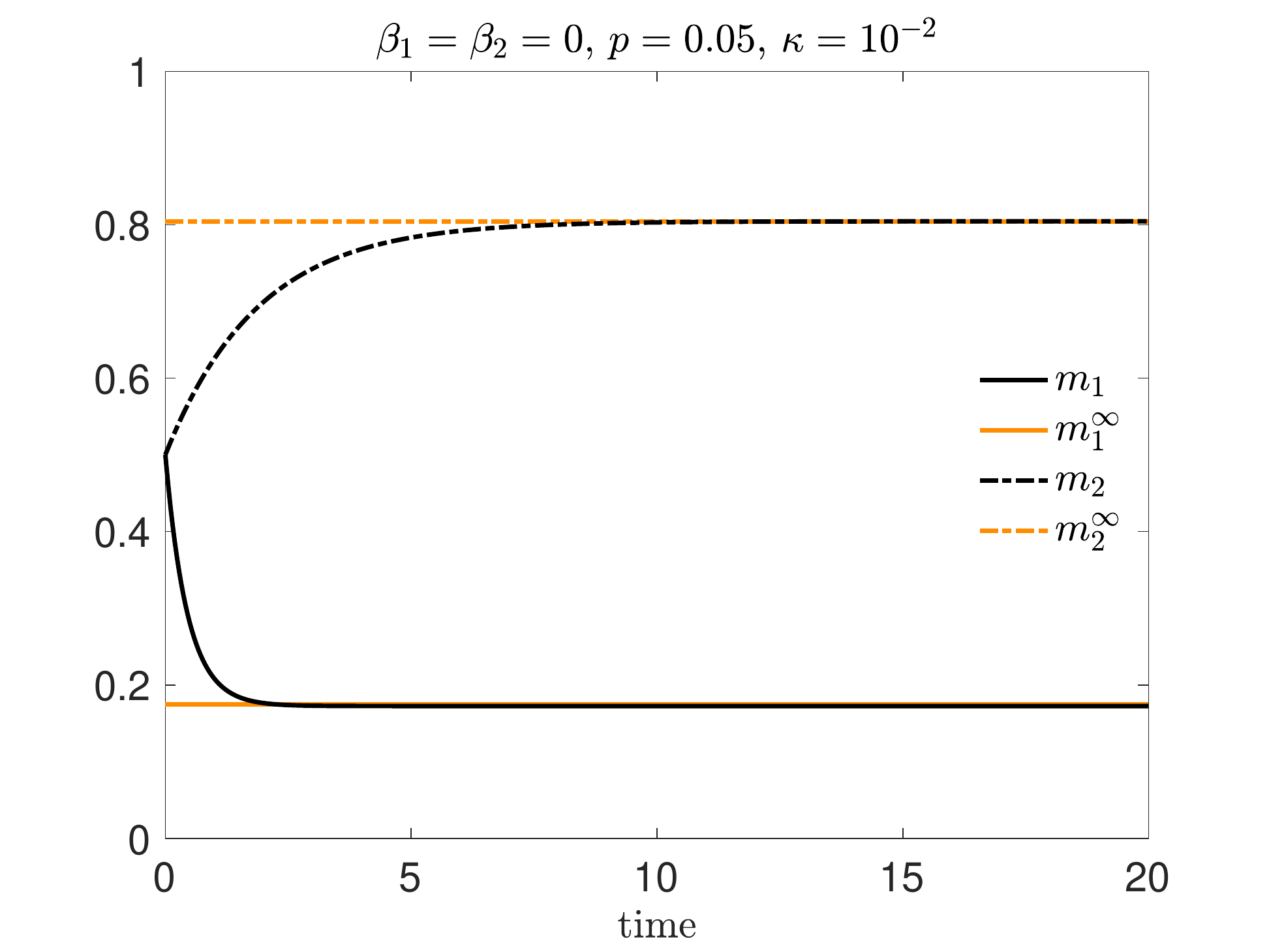} \\
\caption{Evolution of the mean speed under the scaling \eqref{eq:scal} in the case $\beta_1=\beta_2=0$, corresponding to excluding lane-switching dynamics. We present the case $\kappa_1 = \kappa_2 = \kappa$ and from left to right we considered the following penalization constants $\kappa = +\infty$, $\kappa = 10^{-1}$ and $\kappa = 10^{-2}$ and a penetration rate $p = 0.05$. We fixed also $\bar v_1 = (1-\rho_1)$, $\bar v_2 = (1-\rho_2)$ and $\rho_1(0) = 1- \rho_2(0) = 0.8$ in \eqref{eq:system_rho} whereas $m_1(0) = m_2(0) = 0.5$. For the acceleration probability $P$  we considered  the constant $\mu = 2$. The solution of \eqref{eq:m12_general} is compared with the analytical large time mean speed $m_i^\infty$, $i = 1,2$, obtained in \eqref{eq:m12_infty}.}
\label{fig:m12_noswitch}
\end{figure}

\begin{figure}
\includegraphics[scale = 0.25]{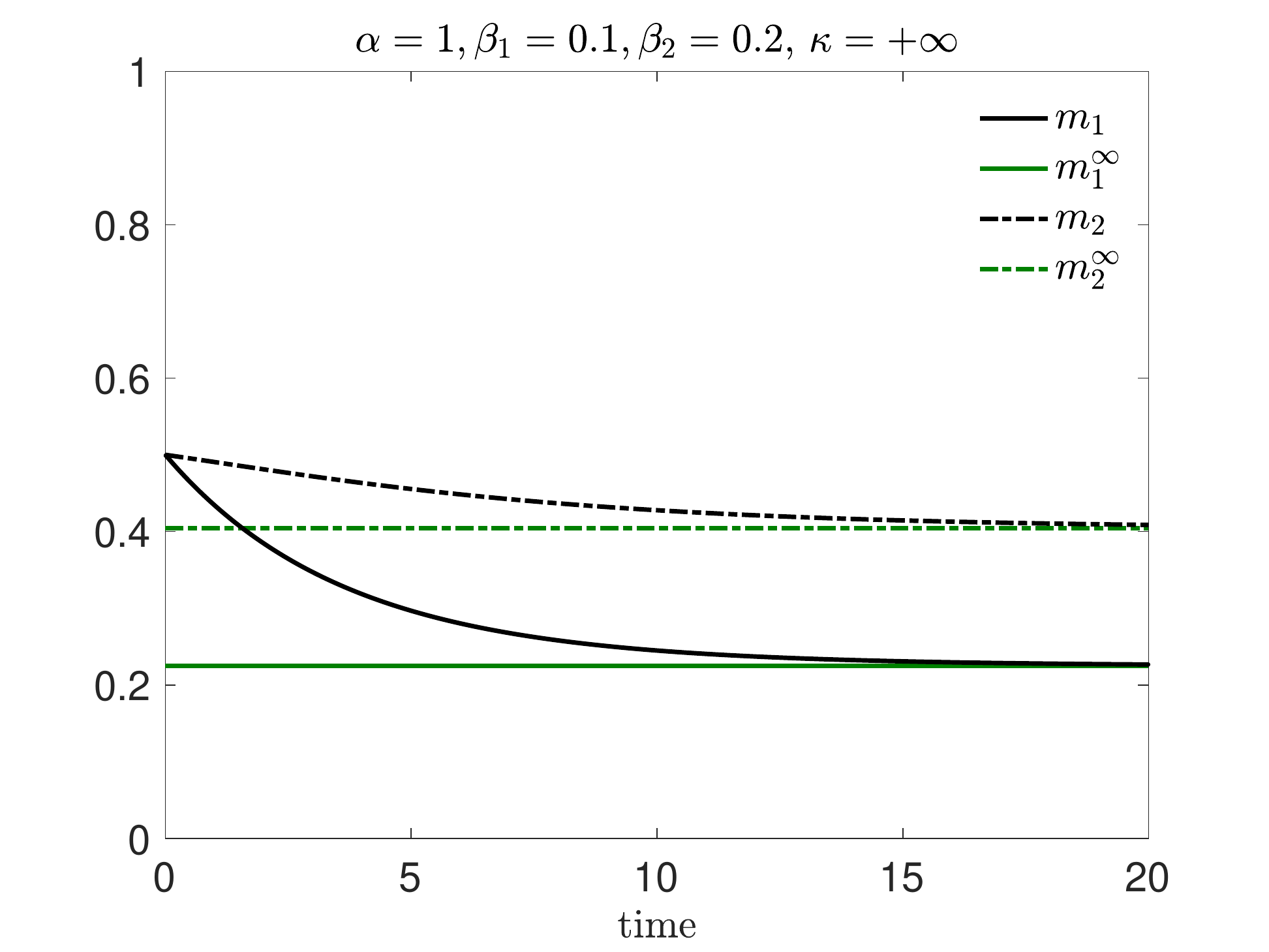} 
\includegraphics[scale = 0.25]{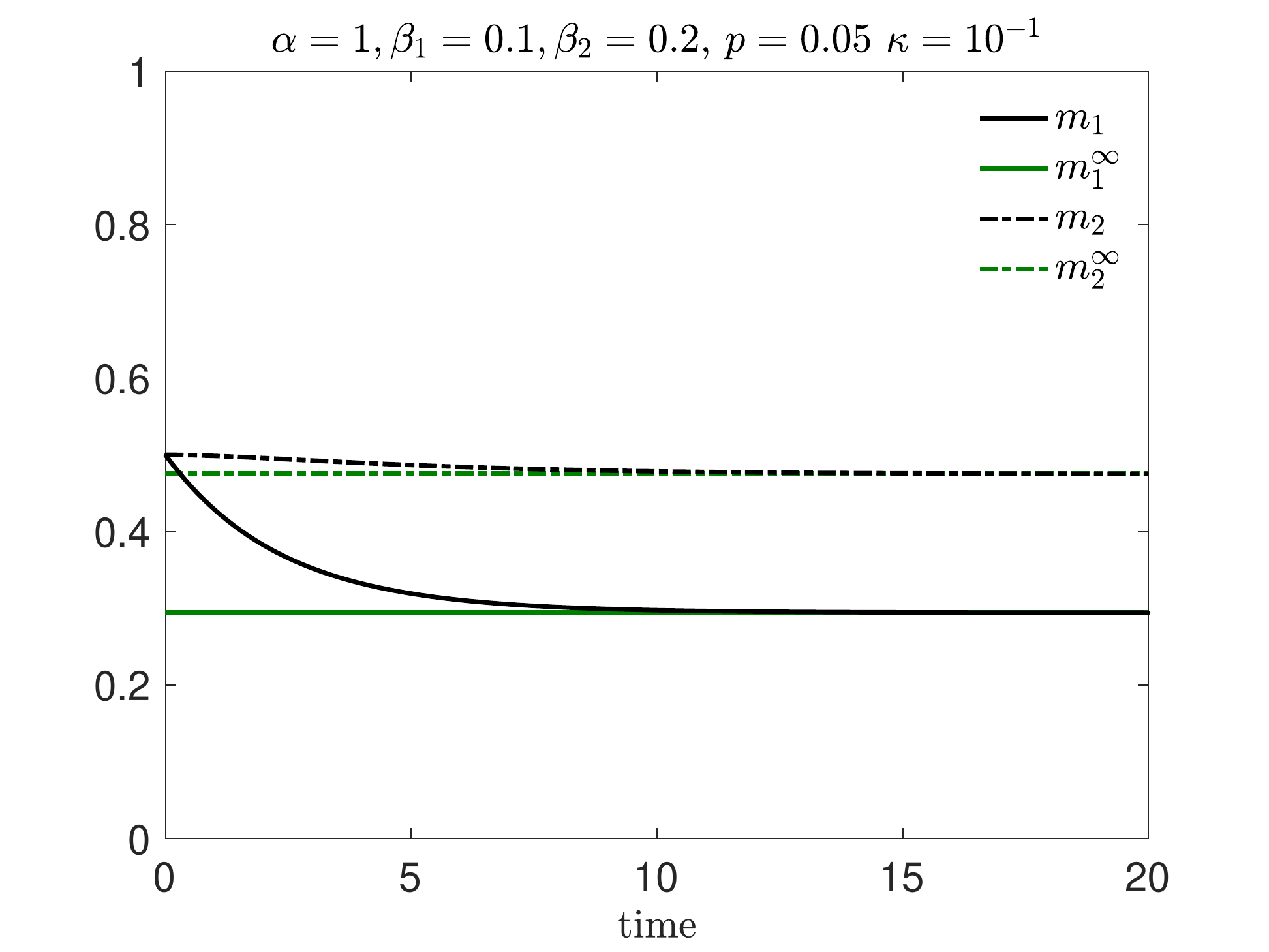} 
\includegraphics[scale = 0.25]{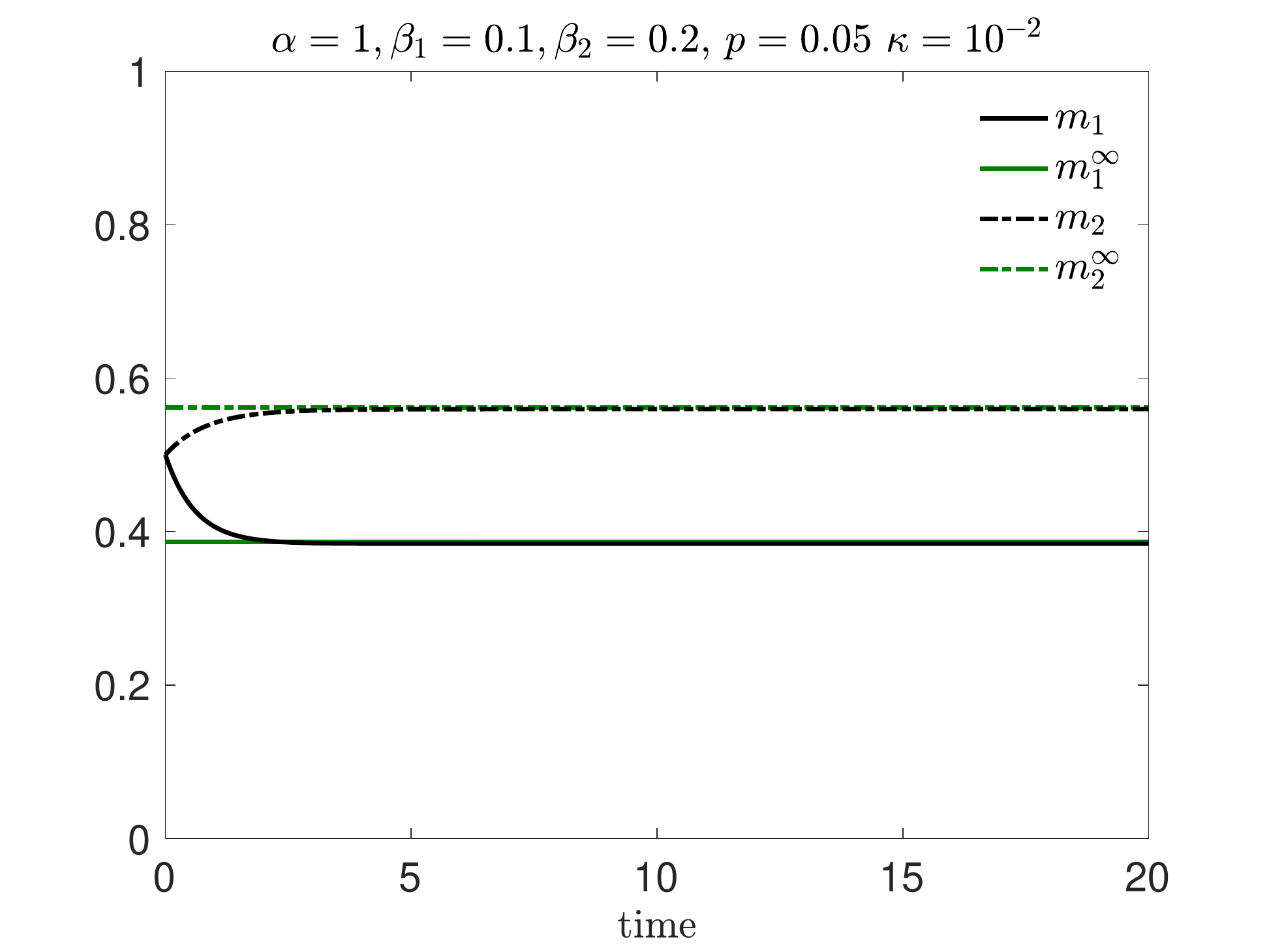} \\
\includegraphics[scale = 0.25]{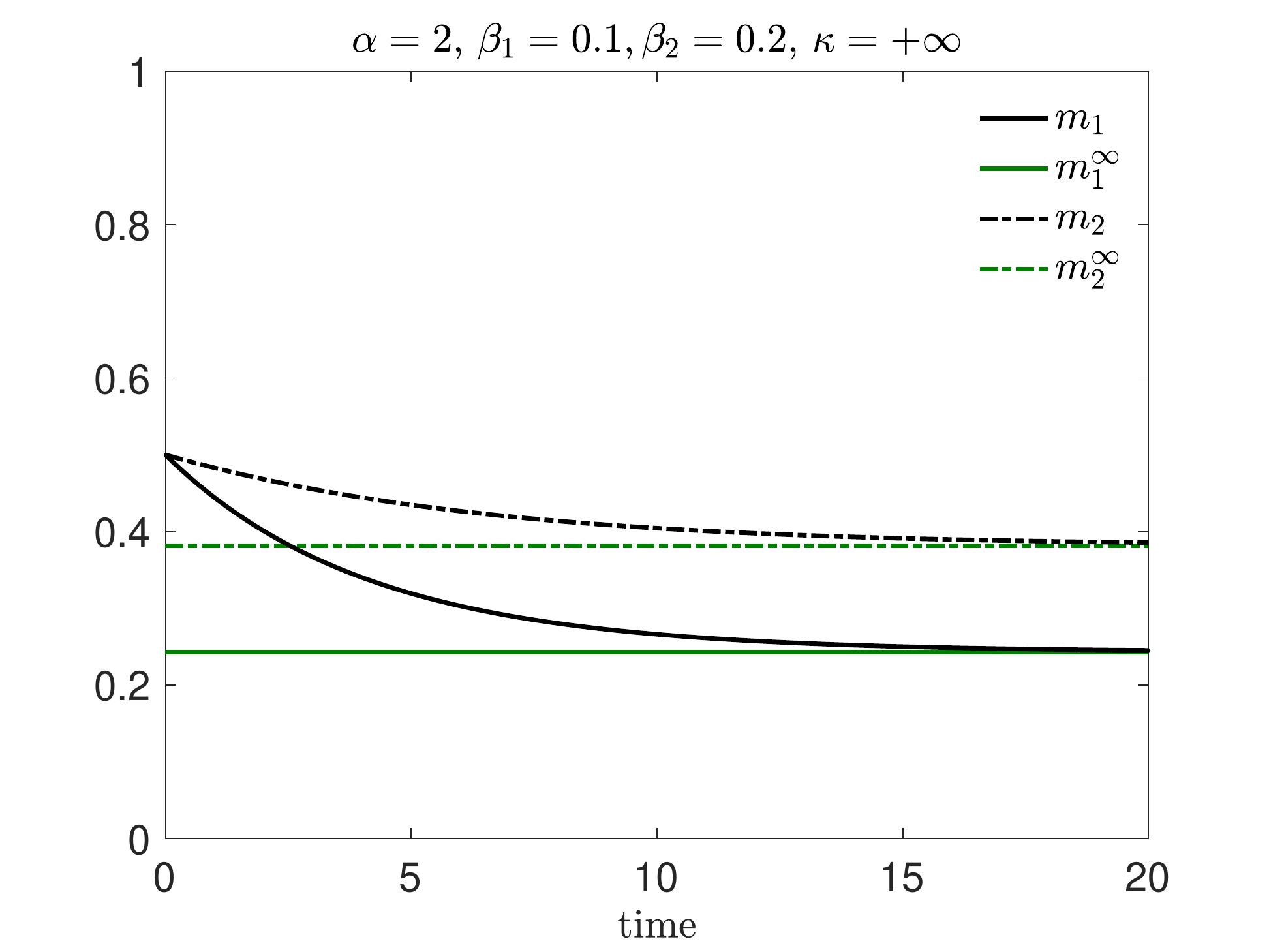}
\includegraphics[scale = 0.25]{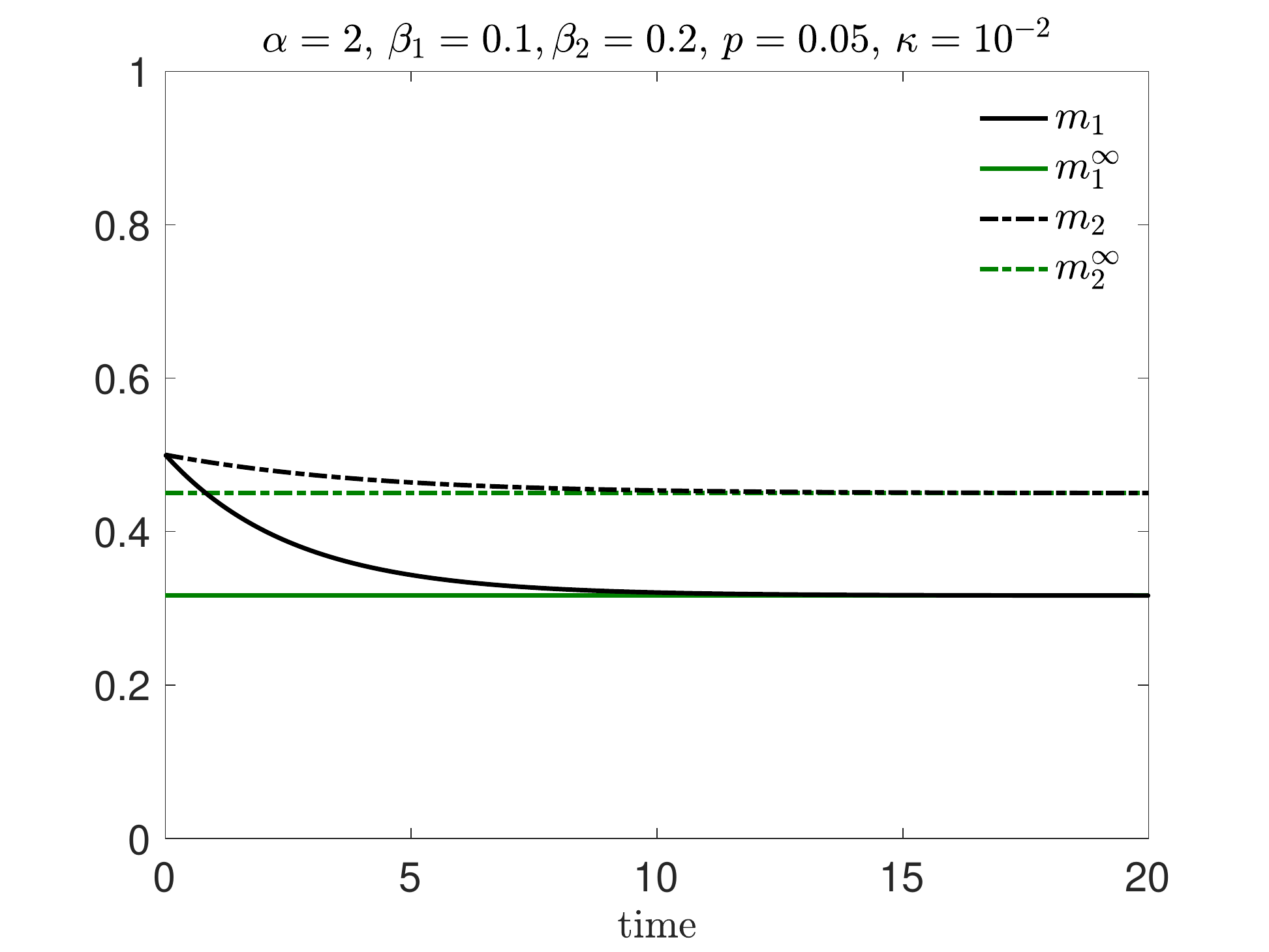}
\includegraphics[scale = 0.25]{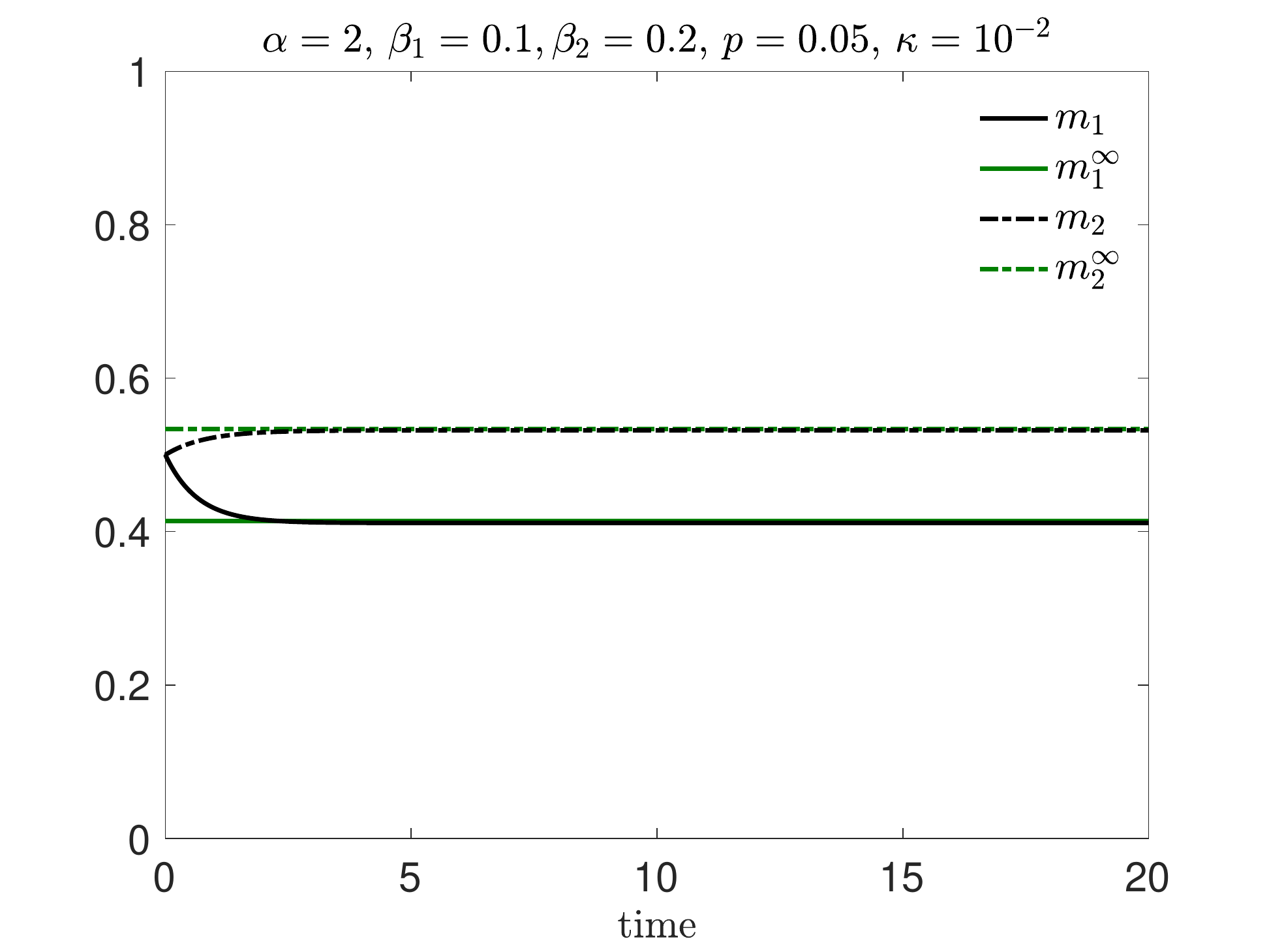}
\caption{Evolution of the mean speed under the scaling \eqref{eq:scal} in the case $\beta_1=0.1$, $\beta_2=0.2$ and $\alpha = 1$ (top row), $\alpha = 2$ (bottom row). We present the case $\kappa_1 = \kappa_2 = \kappa$ and from left to right we considered the following penalization constants $\kappa = +\infty$, $\kappa = 10^{-1}$ and $\kappa = 10^{-2}$ and a penetration rate $p = 0.05$. We fixed also $\bar v_1 = (1-\rho_1)$, $\bar v_2 = (1-\rho_2)$, initial densities $\rho_1(0) = 1- \rho_2(0) = 0.8$, and the initial mean speeds  $m_1(0) = m_2(0) = 0.5$. For the acceleration probability $P(\rho_i)$ we considered $\mu = 2$. We compared the solution of \eqref{eq:m12_general} with the analytical large time $m_i^\infty$, $i = 1,2$ in \eqref{eq:m12_infty}.}
\label{fig:m12_switch}
\end{figure}

\begin{figure}
\centering
\includegraphics[scale = 0.35]{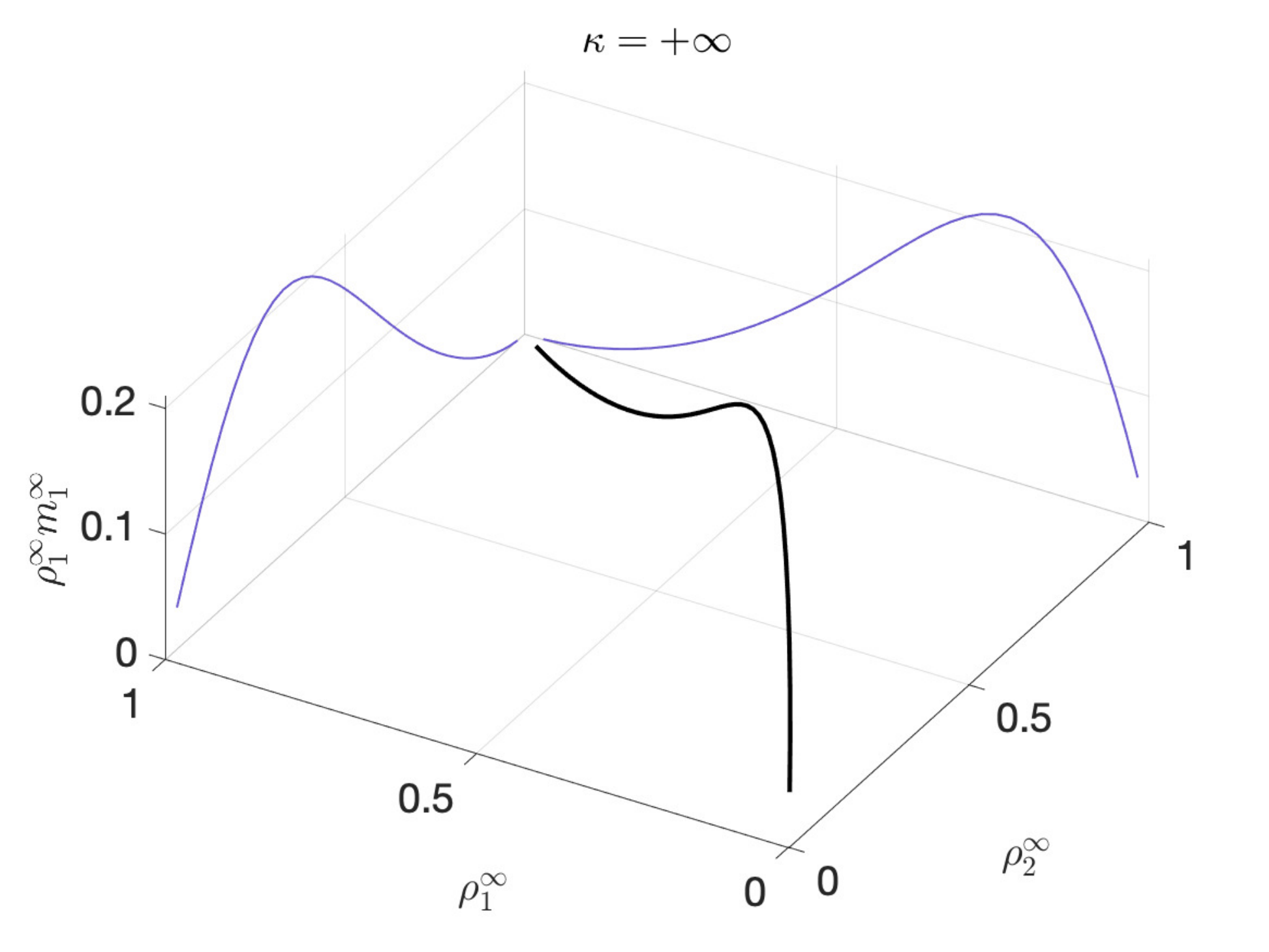}
\includegraphics[scale = 0.35]{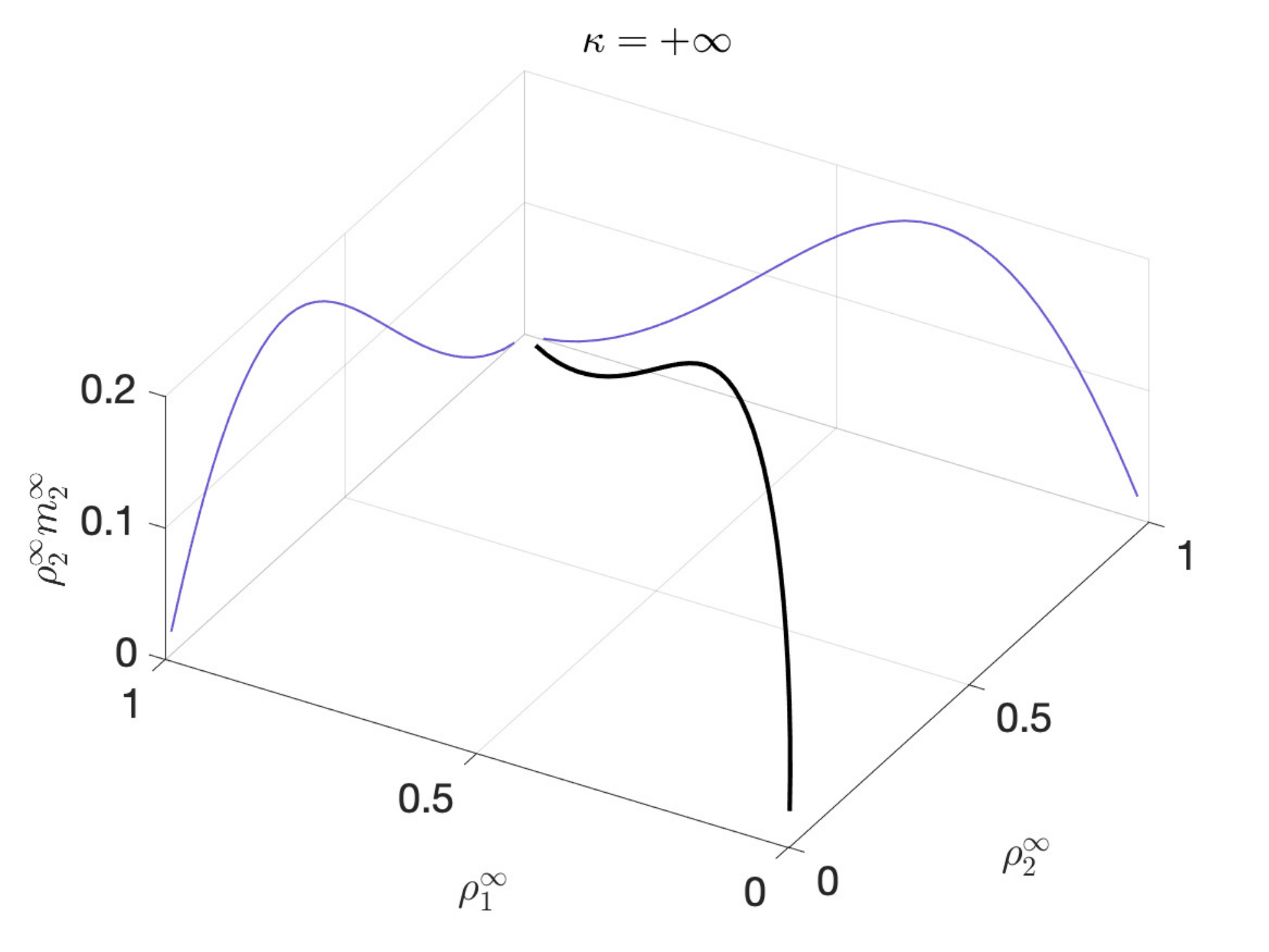} \\
\includegraphics[scale = 0.35]{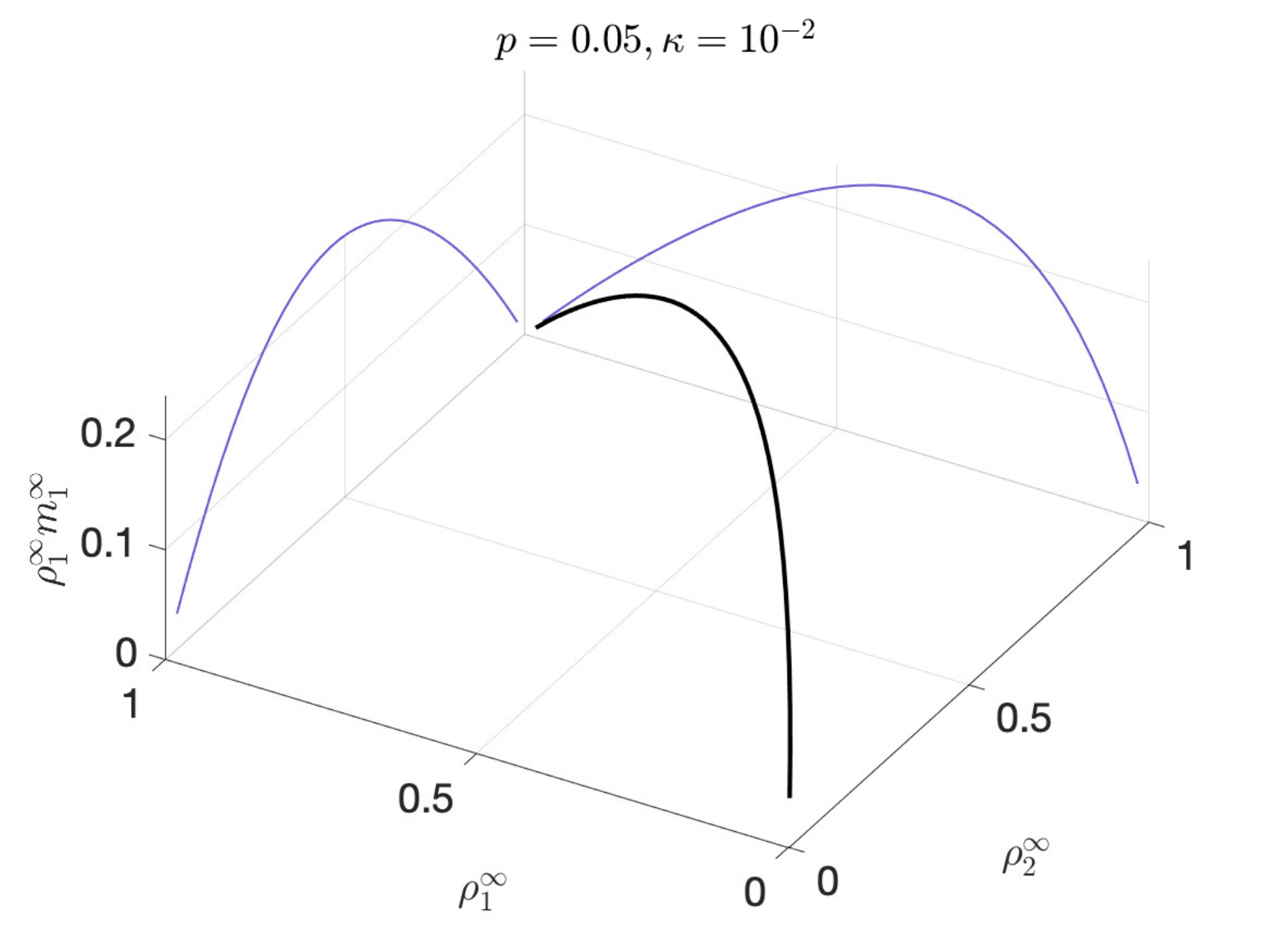}
\includegraphics[scale = 0.35]{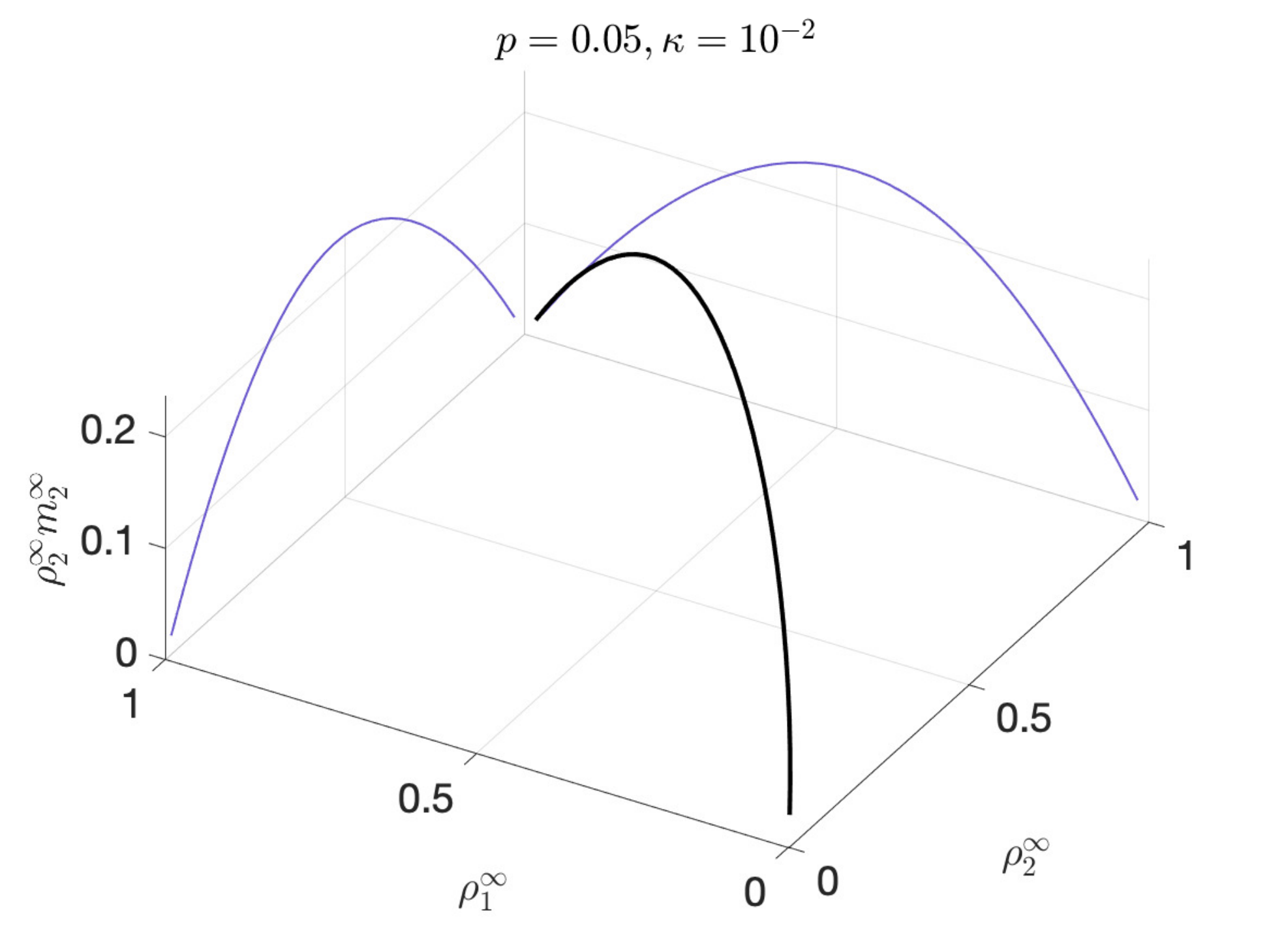} 
\caption{Traffic diagrams of the first lane (left column) and of the second lane (right column) obtained from \eqref{eq:m12_infty} and $\alpha = 2$. Furthermore we present the unconstrained case $ \kappa_1 = \kappa_2 = \kappa = +\infty$ (top row) and the constrained case with $\kappa_1 = \kappa_2 = \kappa = 10^{-2}$ (bottom row). In both cases we considered $\beta_1 = 0.1$, $\beta_2 = 0.2$}
\label{fig:traffic_diag}
\end{figure}

\begin{figure}
\centering
\includegraphics[scale = 0.35]{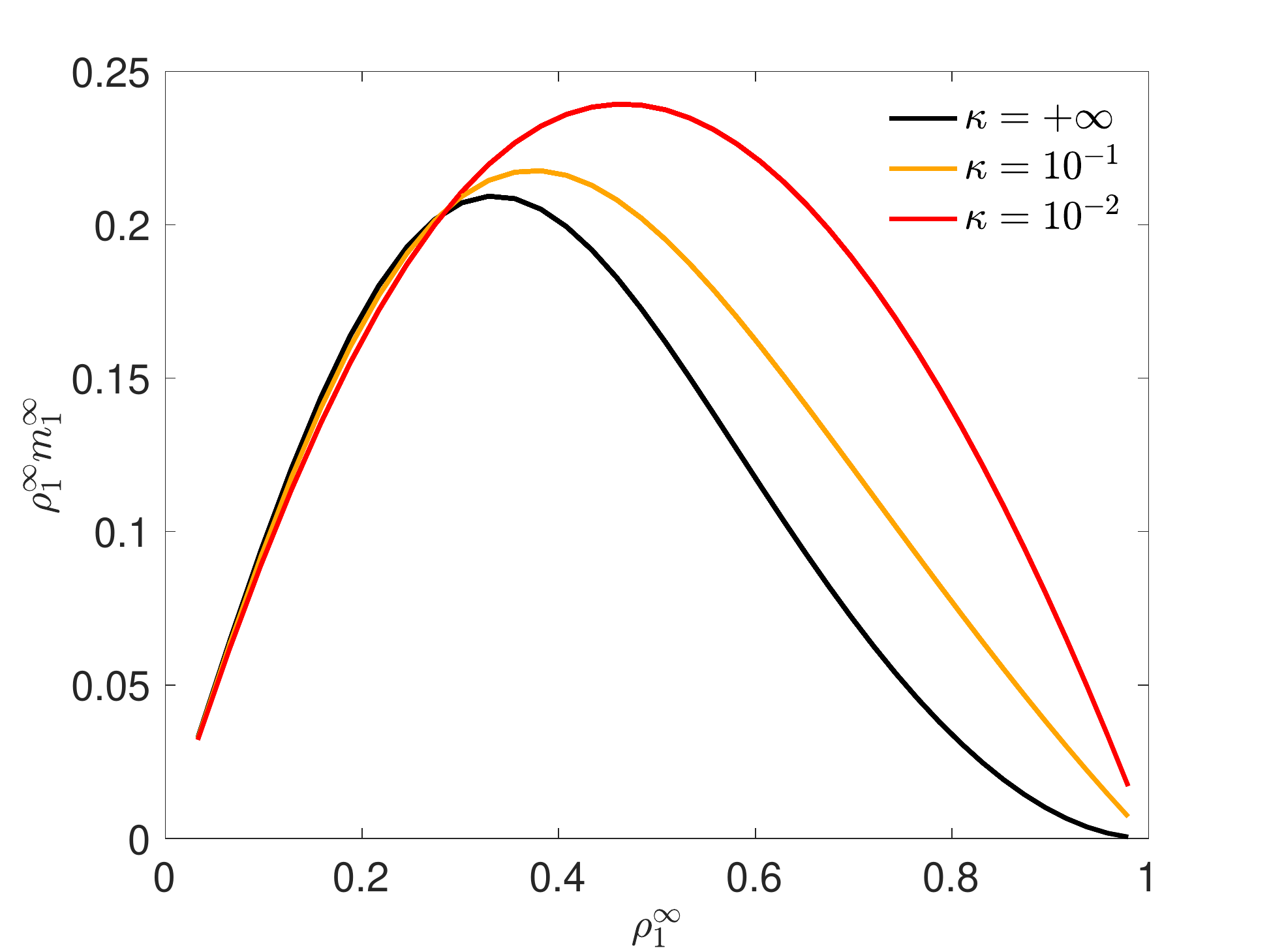}
\includegraphics[scale = 0.35]{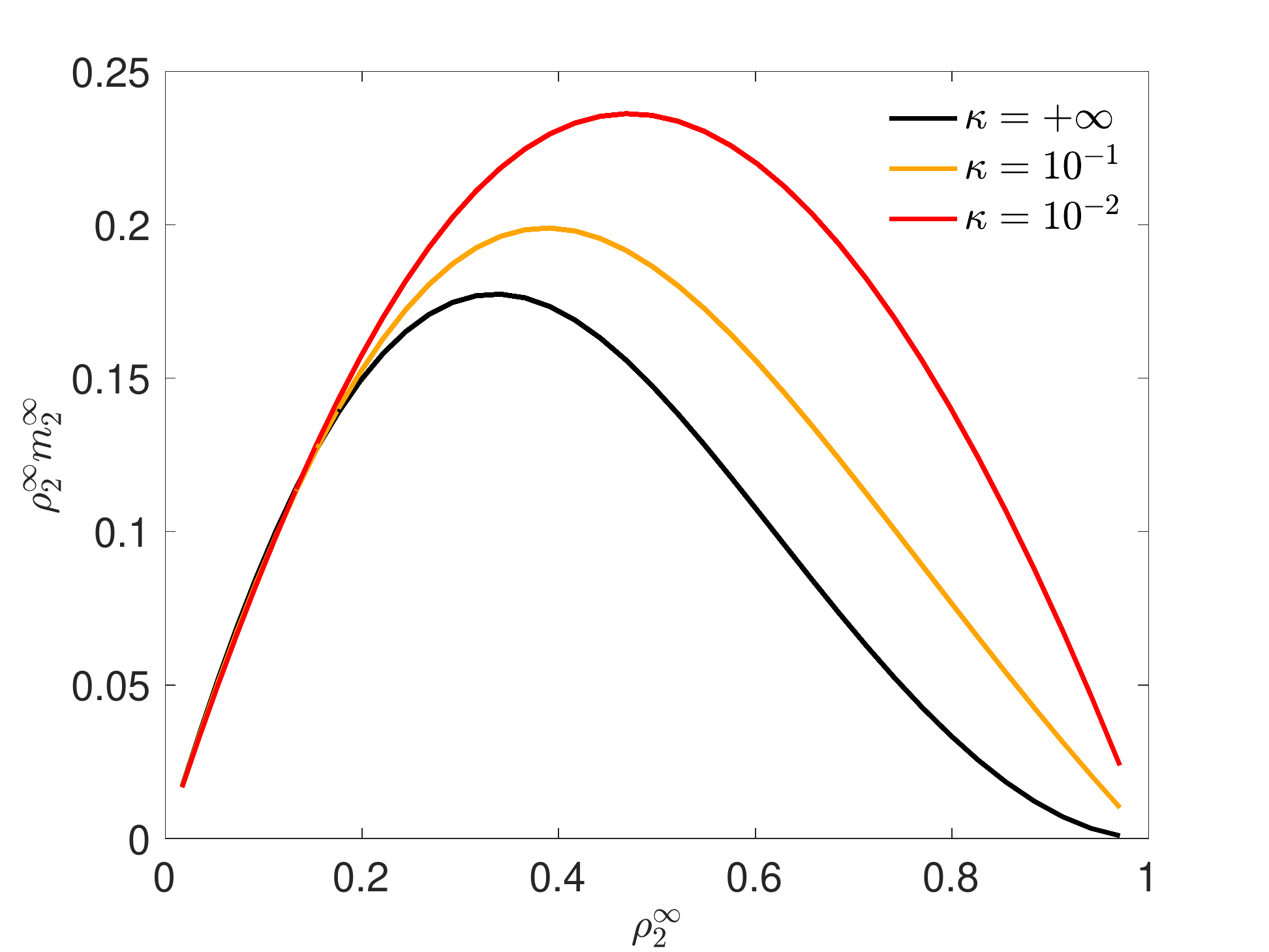} 
\caption{Traffic diagrams of the first lane (left column) and of the second lane (right column) for penetration rate $p = 0.05$ and different penalizations $\kappa = 10^{-2},10^{-1}$ and parameters conformal to Figure \ref{fig:traffic_diag}. The unconstrained case corresponds to $\kappa = +\infty$.}
\label{fig:traffic_diag2} 
\end{figure}

The evolution of the mean speeds in each lane is therefore coupled, in principle, to the evolution of the lanes' densities. In order to get an analytic insight on the asymptotic mean speeds, since \eqref{eq:system_rho} is decoupled from the evolution of the mean speeds, we can suppose that $\rho_1(0)$, $\rho_2(0)$ are at equilibrium such that $\rho_1(0) = \rho^\infty_1>0$, $\rho_2(0)= \rho^\infty_2>0$ for all $t\ge 0$. To simplify notations in the following we denote by $\rho_1,\rho_2$ the equilibrium densities.  Furthermore, we introduce the new time-scale $\tau = t/\gamma$ and we introduce the following scaling
\begin{equation}\label{eq:scal}
\beta_i\rightarrow \beta_i\gamma,\qquad \nu_i = \kappa_i\gamma,\qquad i = 1,2
\end{equation}
to obtain in the limit $\gamma\rightarrow 0^+$ the evolution of the rescaled mean speeds $m_1(t/\gamma)$, $m_2(t/\gamma)$ given by
\begin{equation}\label{eq:m12}
\begin{cases}\vspace{0.25cm}
\dfrac{d}{d\tau} \rho_1 m_1(\tau)=& \dfrac{\rho_1^2}{2} \left[P(\rho_1)  - \left(P(\rho_1) + (1-P(\rho_1))^2\right)m_1 + p_1^* (\bar v_1 - m_1)\right]  \\
&- \beta_1(1-\rho_2)^\alpha \rho_1 m_1 + \beta_2(1-\rho_1)^\alpha \rho_2m_2  \\
\dfrac{d}{d\tau} \rho_2 m_2(\tau) =&\dfrac{\rho_2^2}{2} \left[P(\rho_2)  - \left(P(\rho_2) + (1-P(\rho_2))^2\right)m_2 + p_2^* (\bar v_2 - m_2)\right]  \\
&- \beta_2(1-\rho_1)^\alpha \rho_2m_2 + \beta_1(1-\rho_2)^\alpha \rho_1m_1  \\
\end{cases}
\end{equation}
where $p^*_i=\dfrac{p}{\kappa_i}$ is the effective penetration rate of the traffic flow \cite{TZ1}. Under the introduced assumptions, the large time mean speeds are therefore given by 
\begin{equation}\label{eq:m12_infty}
\begin{split}\vspace{0.25cm}
m_1^\infty &= \dfrac{1}{A_1A_2C_1C_2 - 1} \left( A_2B_1C_1C_2 +C_2B_2 \right)\\
m_2^\infty &=\dfrac{1}{A_1A_2C_1C_2 - 1} \left( A_1B_2C_1C_2 +C_1B_1 \right)
\end{split}
\end{equation}
where
\[
A_i = P(\rho_i) + (1-P(\rho_i))^2 + p^*_i + \dfrac{2\beta_i}{\rho_i} (1-\rho_j)^\alpha, \qquad B_i = P(\rho_i) + p^*_i \bar v_i,\qquad
C_i = \dfrac{\rho_i^2}{2\beta_j \rho_j (1-\rho_i)^\alpha}
\]
with $i,j = 1,2$, $i\ne j$. 
In Figure \ref{fig:m12_noswitch} and Figure \ref{fig:m12_switch} we test the consistency of the large time mean speeds $m_i^\infty$ in \eqref{eq:m12_infty} with the solution of \eqref{eq:m12_general} under the scaling \eqref{eq:scal} with $\gamma = \Delta t = 10^{-3}$. To obtain asymptotic densities we first solved \eqref{eq:system_rho} supposing $\rho_1(0) = 1-\rho_2 (0)= 0.8$. 

From \eqref{eq:m12_infty} one observes that for infinite penalization $\kappa_i\rightarrow +\infty$ no effect of the control is present in the $i$th lane since $p^*_i=0$. This corresponds to the case of zero penetration rate $p=0$, i.e. in the absence of driver-assist vehicles in the stream. On the other hand, in case of negligible penalization, we obtain for the $i$th lane 
\[
m_i^\infty(\rho_i) = \bar v_i, 
\]
even if $\kappa_j\gg1$, $j\ne i$, meaning that we are capable to obtain lane alignment of velocities also if  lane switching dynamics are allowed. 

In the case $\beta_1 = \beta_2 =0$, the densities $\rho_1$ and $\rho_2$ become conserved quantities and \eqref{eq:m12} reduces to
\[\vspace{0.25cm}
\begin{cases}
\dfrac{d}{d\tau} m_1 = P(\rho_1) - (P(\rho_1)+(1-P(\rho_1))^2)m_1 + p_1^*(\bar v_1 - m_1),\\
\dfrac{d}{d\tau} m_2 = P(\rho_2) - (P(\rho_2)+(1-P(\rho_2))^2)m_2 + p_2^*(\bar v_2 - m_2),
\end{cases}
\]
and therefore the large time mean velocities are given by 
\begin{equation}
\label{eq:mtilde12}
\begin{split}
\tilde m_1^\infty(\rho_1) = \dfrac{P(\rho_1) + p_1^* \bar v_1}{P(\rho_1) + (1-P(\rho_1))^2 +p_1^*}, \qquad \tilde m_2^\infty(\rho_2) = \dfrac{P(\rho_2) + p_2^*\bar v_2}{P(\rho_2) + (1-P(\rho_2))^2 + p_2^*}.
\end{split}
\end{equation}
This solution is coherent with the single lane mean speed of \cite{TZ1}.

In order to get a deeper insight on the effects of lane switching dynamics we consider in Figure \ref{fig:traffic_diag}-\ref{fig:traffic_diag2} the flux-density diagrams both in the unconstrained and constrained regimes with $\alpha = 2$. Those diagrams establish a relationship between macroscopic flux of vehicles to the traffic density in homogeneous conditions. Empirical observation, that are common to all measured traffic diagrams, suggest, that the mean speed is nearly constant  in the free flow regime corresponding to small densities and that the fluxes decrease once the density is higher than a critical density. 

To obtain the fundamental diagrams,  we solve \eqref{eq:system_rho}  with initial densities in $[0,1]$, such that $\rho_1(0)+ \rho_2(0) = \rho \in [0,1]$ and determine the associated asymptotic densities. Then, \eqref{eq:m12_infty} give the large time mean speed and the lane's flux as $\rho_i^\infty m_i^\infty$, $i =1,2$. 

In the second row of Figure \ref{fig:traffic_diag} we plot the resulting traffic diagram in presence of vehicles with the introduced driver-assist technology. As previously discussed, we are indeed capable to drive the macroscopic flux also in presence of switching dynamics to induce alignment of velocities towards a recommended speed, which can also be based on the external information on the local traffic density. The action of the introduced controls is able to increase the flux as we can observe in Figure \ref{fig:traffic_diag2}.

The evolution of the energy of the first and second lane
\[
E_1(t) = \dfrac{1}{\rho_1}\int_0^1 v^2 f_1(v,t)dv,\qquad E_2(t) = \dfrac{1}{\rho_2}\int_0^1 v^2f_2(v,t)dv
\]
are obtained from \eqref{eq:1_weak}-\eqref{eq:2_weak} setting $\varphi(v)= v^2$. Their evolution is coupled to the evolution of both the lanes' densities and mean speeds therefore, to have an analytic insight on their behavior we consider $D \equiv 0$ in the binary interaction term \eqref{eq:v1-2_control}, we consider $\rho_i(0) = \rho_i^\infty$, $m_i(0) = m_i^\infty$ and we introduce the scaling \eqref{eq:scal} to obtain in the limit $\gamma\rightarrow 0^+$ the evolution of the rescaled energies $E_1(t/\gamma)$, $E_2(t/\gamma)$ given by 
\[
\begin{cases}
\dfrac{d}{d\tau} \rho_1 E_1(\tau) =& \rho_1^2 \left[ P(\rho_1) (m_1-E_1) + (1-P(\rho_1))^2 (m_1^2 - E_1^2) + p^*_1(\bar v_1m_1 - E_1) \right] \\
& - \beta_1 (1-\rho_2)^\alpha \rho_1 E_1 + \beta_2 (1-\rho_1)^\alpha \rho_2E_2 \\
\dfrac{d}{d\tau} \rho_2 E_2(\tau) =& \rho_2^2 \left[ P(\rho_2) (m_2- E_2) + (1-P(\rho_2))^2 (m_2^2 - E_2^2) + p^*_2(\bar v_2m_2 - E_2) \right] \\
& - \beta_2 (1-\rho_1)^\alpha \rho_2 E_2 + \beta_1 (1-\rho_2)^\alpha \rho_1 E_1,  \\
\end{cases}
\]
with $\tau = t/\gamma$. We observe that for $\kappa_i \rightarrow 0^+$ we have
\[
E_i  \rightarrow \bar v_i^2, 
\]
and therefore $E_i - m_i^2 \rightarrow 0$ such that the large time solution of the space homogeneous version of \eqref{eq:1_weak}-\eqref{eq:2_weak} are Dirac delta $\delta(v-\bar v_i)$.

\section{Hydrodynamic regime}\label{sec:hydro}

\subsection{Homogeneous equilibria of automated traffic}\label{sec:hom}

Before discussing hydrodynamic limits of the introduced traffic dynamics we study the asymptotic distribution of speed and variance of the homgeneous equation in more detail and give an approximation for the stationary solution of the homogeneous equation.  We consider  equation (\ref{eq:1_weak}) and (\ref{eq:2_weak}) with $\beta_1  = \beta_2 = 0$,  i.e.
\begin{equation}
\label{eq:hom_1}
\dfrac{d}{dt} \int_0^1\varphi(v) f_i(v,t)dv = \mathbb E_{\Theta,\eta_i} \left[ \int_0^1 \int_0^1 \left( \varphi(v^\prime) - \varphi(v) \right) f_i(v,t)f_i(w,t)\,dv\,dw \right], \qquad i = 1,2. 
\end{equation}
We introduce the scaled homogeneous distribution function $f(v,\tau) = f(2\tau/\gamma,v)$ and consider the large time limit $\gamma\rightarrow 0^+$. 
Although the large time density and mean speeds could be determined in the previous section, the exact shape  of the stationary distribution function is not available. 
 A reduction of complexity can be achieved in the quasi-invariant interaction limit $\gamma, \sigma^2_i \rightarrow 0^+$, for $i=1,2$. This limit is reminiscent of the grazing collision limit of gas dynamics, see \cite{FPTT,PTZ,Vill}. 

On the new time scale we  study 
\begin{equation}
\label{eq:hom_2}
\dfrac{d}{d\tau} \int_0^1 \varphi(v)f_i(v,\tau)dv = \dfrac{1}{\gamma}\mathbb E_{\Theta,\eta_i} \left[ \int_0^1 \int_0^1 (\varphi(v^\prime)-\varphi(v))f_i(v,\tau)f_i(w,\tau)\,dv \, dw\right], \qquad i = 1,2. 
\end{equation}
Assuming a sufficient regularity of $\varphi\in C^3([0,1])$ we can perform a Taylor expansion
\[
\varphi(v^\prime)-\varphi(v) = \varphi^\prime(v)(v^\prime-v) + \dfrac{1}{2}\varphi^{\prime\prime}(v)(v^\prime-v)^2 + \dfrac{1}{6}\varphi^{\prime\prime\prime}(\tilde v)(v^\prime-v)^3,
\]
with $\tilde v \in (\min\{v,v^\prime\},\max\{v,v^\prime\})$. Now we can plug the expansion in \eqref{eq:hom_2} and use the scaling 
\begin{equation}\label{eq:quasi-invariant}
\qquad \nu_i=\kappa_i \gamma, \qquad \sigma^2_i = \gamma \lambda_i.
\end{equation}
Then, $\gamma \rightarrow 0^+$ gives 
\[
\begin{split}
\dfrac{d}{d\tau} \int_0^1 \varphi(v)f_i(v,\tau)dv =& \int_0^1 \int_0^1   \varphi^\prime(v)\left[I(v,w,\rho_i) + \dfrac{p}{\kappa_i}(\bar v_i(\rho_i) - v) \right]f_i(v,\tau)f_i(w,\tau)dv\,dw \\
&+ \dfrac{\lambda_i \rho_i}{2} \int_0^1 \varphi^{\prime\prime}(v)D^2(v,\rho_i)f_i(v,\tau)dv + R_{\varphi}(f_i,f_i). 
\end{split}
\]
It is possible to prove that if the random variables $\eta_i$ are sufficiently regular, namely with bounded third order, i.e. $\mathbb E_{\eta_i}[|\eta|^3]<+\infty$, the remainder term $R_{\varphi} \rightarrow 0$ for $\gamma\rightarrow 0^+$. A detailed analysis can be found in \cite{CPT,T,TZ2}. 

This gives a Fokker-Planck equation  with nonlocal drift and nonconstant diffusion
\[
\partial_\tau f_i = \partial_v \left[ \left(\dfrac{p\rho_i}{\kappa_i}v- \int_0^1  \left(I(v,w,\rho_i) + \dfrac{p}{\kappa_i}\bar v_i \right) f_i(w,\tau)dw \right)  f_i(v,\tau) + \dfrac{\lambda \rho_i}{2}\partial_v \left( D^2(v,\rho_i) f_i \right)\right].
\]
Compare \cite{IKM} for another  Vlasov-Fokker-Planck approach to multilane traffic.
Asymptotic states can be explicitly computed for choices of the diffusion function $D(v,\rho_i)$. If we chose 
\[
D(v,\rho_i) = a(\rho_i)\sqrt{v(1-v)}, \qquad a(\rho_i)\ge 0,
\]
the stationary distribution $f_i^{\infty}(v)$ are such that
\begin{equation}
\label{eq:condition_inf}
\dfrac{\lambda_i}{2} \partial_v^2 (D^2(v,\rho_i) f_i^{\infty}) - (1+\dfrac{p}{\kappa_i})\partial_v \left( (m^\infty_i - v)f_i^\infty \right) = 0,
\end{equation}
where $m^\infty_i$ are given in  \eqref{eq:mtilde12}. In particular, the solution of \eqref{eq:condition_inf} reads
\begin{equation}
\label{eq:f_infty}
f_i^\infty(v) = \dfrac{v^{\frac{2(1+p/\kappa_i)}{\lambda_i a(\rho_i)}m_i^\infty(\rho_i) - 1}(1-v)^{\frac{2(1+p/\kappa_i)}{\lambda_i a(\rho_i)}(1-m_i^\infty(\rho_i)) - 1}}{\mathbb{B}\left( \frac{2(1+p/\kappa_i)}{\lambda_i a(\rho_i)}m_i^\infty(\rho_i), \frac{2(1+p/\kappa_i)}{\lambda_i a(\rho_i)}(1-m_i^\infty(\rho_i)) \right)},
\end{equation}
where $\mathbb B(x,y) = \int_0^1 t^{x-1}(1-t)^{y-1}dt$ is the beta function. We observe that  \eqref{eq:f_infty} are beta probability distribution $\mathbb {B}(I_i,J_i)$ with
\[
I_i = \dfrac{2(1+\frac{p}{\kappa_i}) m_i^\infty(\rho_i)}{\lambda_i a(\rho_i)}, \qquad J_i = \dfrac{2(1+\frac{p}{\kappa_i})(1-m_i^\infty(\rho_i))}{\lambda_i a(\rho_i)},  
\] 
and it is easily checked that if the random variable $X\sim f_i^\infty$ describes the asymptotic distribution of vehicles' speeds we have
\begin{equation}\begin{split}
\mathbb E[X] &= \dfrac{I_i}{I_i + J_i} = m_i^\infty(\rho_i), \\
\textrm{Var}(X) &= \dfrac{I_i J_i}{(I_i+J_i)^2 (I_i+J_i+1)} = \dfrac{\lambda_i a^2(\rho_i)}{2+\lambda_i a^2(\rho_i) + 2\frac{p}{\kappa_i}} m_i^\infty(\rho_i)(1-m_i^\infty(\rho_i)). 
\end{split}
\end{equation}
Remarkably enough, beta-type velocity distributions are often observed in experimental literature. We point the interested reader to \cite{HTVZ2,NHJ} for extensive discussions on a data-driven approach based on experimental data.

\subsection{First order hydrodynamics}
We may recover transport models from the inhomogeneous system \eqref{eq:multi} by introducing a hyperbolic scaling of space and time
\[
t\rightarrow \dfrac{t}{\epsilon},\qquad x \rightarrow \dfrac{x}{\epsilon} \qquad (0<\epsilon\ll 1),
\]
to define proper time and space observable at the fluid regime. In the following we will denote with $f(x,v,t)$ the scaled distribution $f(x/\epsilon,v,t/\epsilon)$. The scaled model reads
\begin{equation}
\label{eq:multi_s1}
\begin{split}
\partial_t \int_0^1 \varphi(v)f_1(x,v,t)dv + \partial_x \int_0^1 v\varphi(v)f_1(x,v,t)dv = \dfrac{1}{\epsilon} \int_0^1 \varphi(v)Q(f_1,f_1)(x,v,t)dv \\
- \dfrac{\beta_1}{\epsilon}(1-\rho_2)^{\alpha}\int_0^1 \varphi(v)f_1(x,v,t)dv + \dfrac{\beta_2}{\epsilon}(1-\rho_1)^\alpha \int_0^1 \varphi(v)f_2(x,v,t)dv,
\end{split}\end{equation}
and
\begin{equation}
\label{eq:multi_s2}
\begin{split}
\partial_t \int_0^1 \varphi(v)f_2(x,v,t)dv + \partial_x \int_0^1 v\varphi(v)f_2(x,v,t)dv = \dfrac{1}{\epsilon} \int_0^1 \varphi(v)Q(f_2,f_2)(x,v,t)dv \\
- \dfrac{\beta_2}{\epsilon}(1-\rho_1)^{\alpha}\int_0^1 \varphi(v)f_2(x,v,t)dv + \dfrac{\beta_1}{\epsilon}(1-\rho_2)^\alpha \int_0^1 \varphi(v)f_1(x,v,t)dv.
\end{split}\end{equation}

We will consider  three different regimes generating several fluid models for the multilane constrained traffic dynamics:
\begin{enumerate}
\item[$i)$] Collision dominated regime: $\beta_1,\beta_2 = o(\epsilon)$
\item[$ii)$] Fast switching regime: $\beta_1,\beta_2 = O(1)$
\item[$iii)$] Slow switching regime: $\beta_1,\beta_2 = O(\epsilon) $.
\end{enumerate}

\subsubsection{Collision dominated regime}\label{sect:collision_dom}
If $\beta_1/\epsilon \rightarrow 0^+$, meaning that $\beta_1,\beta_2 = o(\epsilon)$, no effect of the lane changing is noticed in the fluid regime. The  dynamics on the fast and on the slow time scale of the scaled two-lane model \eqref{eq:multi_s1}-\eqref{eq:multi_s2} are 

\begin{equation}
\begin{cases}\vspace{0.25cm}
\partial_t \displaystyle\int_0^1 \varphi(v)f_i(x,v,t)dv = \dfrac{1}{\epsilon} \displaystyle\int_0^1 \varphi(v)Q(f_i,f_i)(x,v,t)dv \\
\partial_t \displaystyle\int_0^1 \varphi(v) f_i(x,v,t)dv + \partial_x\displaystyle \int_0^1 v\varphi(v)f_i(x,v,t)dv = 0.
\end{cases}
\end{equation}
If $\epsilon\rightarrow 0^+$, in the quasi-invariant regime \eqref{eq:quasi-invariant} we may take advantage of the space-homogeneous equilibrium solution based on the derived Fokker-Planck operator 
\[
f_i(x,v,t) \approx \rho_i(x,t) f_i^\infty(v), \qquad i = 1,2,
\]
which provides us with the  local equilibrium closure. Hence, going to the transport step we obtain the macroscopic decoupled system of conservation laws with $i=1,2$
\[
\partial_t \left(\rho_i \displaystyle\int_0^1 \varphi(v)f_i^\infty(v)dv \right)+ \partial_x \left(\rho_i \displaystyle\int_0^1 v\varphi(v)f_i^\infty(v)dv \right) = 0 
\]
Since, in this regime, density is conserved, corresponding to the zeroth moment of the kinetic distribution function, a closed hydrodynamic system of equations can be obtained by choosing $\varphi(v) = 1$

\[
\partial_t \rho_i(x,t) + \partial_x \left( \rho_i(x,t) m_1^\infty(\rho_i) \right)=0  \\
\]
being $m_i^\infty(\rho_i) = \int_0^1 v f_i^\infty(v)dv$, $i=1,2$. Since in the regime $\beta_1,\beta_2 = o(\epsilon)$ no effect of the lane changing is detected on the kinetic dynamics, we can take advantage of \eqref{eq:m12_inf}. The explicit system of first order hydrodynamic equations is then derived
\begin{equation}
\partial_t \rho_i(x,t) + \partial_x \left( \rho_i(x,t) \dfrac{P(\rho_i) + \frac{p}{\kappa_i} \bar v_i}{P(\rho_i) + (1-P(\rho_i))^2 + \frac{p}{\kappa_i}} \right)=0 .
\end{equation}

This macroscopic model is coherent with the one proposed in \cite{TZ1} since the frequency of lane changing dynamics is too weak compared with ithe nteraction frequency. The obtained decoupled system of first-order hydrodynamic traffic model has the flux $\mathcal F = (\rho_1 \tilde m_1^\infty,\rho_2 \tilde m_2^\infty)$, where $\tilde m_1^\infty$ has been defined in \eqref{eq:mtilde12}. This flux is not necessarily concave for all $\rho_i \in [0,1]$ as for classical Lighthill-Whitham-Richards models for vehicular traffic. Furthermore the flux $\mathcal F$ depends on the microscopic control strategy defined by the portion of driver-assist vehicles in the traffic. 

\subsubsection{Fast switching regime}\label{sect:slow}

In this regime we consider $\beta_1,\beta_2 = O(1)$ and the scale of interactions is the same as the  one of lane switching. Therefore, in the limit $\epsilon \rightarrow 0^+$ we may split the dynamics  \eqref{eq:multi_s1}-\eqref{eq:multi_s2} as follows
\begin{equation}
\label{eq:1}
\begin{cases}
& \partial_t \displaystyle \int_0^1 \varphi(v) f_1(x,v,t)\,dv = \dfrac{1}{\epsilon}\displaystyle \int_0^1 \varphi(v)Q(f_1,f_1)(x,v,t)dv  \\
& \qquad- \dfrac{\beta_2}{\epsilon} (1-\rho_1 )^\alpha \displaystyle \int_0^1 \varphi(v)f_2(x,v,t) + \dfrac{\beta_1}{\epsilon} (1-\rho_2)^\alpha \displaystyle \int_0^1 \varphi(v)f_1(x,v,t)dv \\
& \partial_t  \displaystyle \int_0^1 \varphi(v)f_1(x,v,t)dv + \partial_x  \displaystyle\int_0^1 v\varphi(v)f_1(x,v,t)dv = 0
\end{cases}
\end{equation}
and 
\begin{equation}
\label{eq:2}
\begin{cases}
& \partial_t \displaystyle \int_0^1 \varphi(v) f_2(x,v,t)\,dv = \dfrac{1}{\epsilon}\displaystyle \int_0^1 \varphi(v)Q(f_2,f_2)(x,v,t)dv  \\
& \qquad- \dfrac{\beta_1}{\epsilon} (1-\rho_2 )^\alpha \displaystyle \int_0^1 \varphi(v)f_1(x,v,t) + \dfrac{\beta_2}{\epsilon} (1-\rho_1)^\alpha \displaystyle \int_0^1 \varphi(v)f_2(x,v,t)dv \\
& \partial_t  \displaystyle \int_0^1 \varphi(v)f_2(x,v,t)dv + \partial_x  \displaystyle\int_0^1 v\varphi(v)f_2(x,v,t)dv = 0. 
\end{cases}
\end{equation}

In the quasi-invariant limit \eqref{eq:quasi-invariant} and $\epsilon \rightarrow 0^+$ this regime makes the lane changing more frequent than interactions between vehicles.
Considering the  terms of order $\frac{1}{\epsilon}$  in the above equations and using $\varphi(v) =1$ we obtain for each $x\in \mathbb R$

\begin{equation}
\label{eq:rho}
\begin{split}
-\beta_2 (1-\rho_1) ^\alpha \rho_2 + \beta_1(1-\rho_2 )^\alpha \rho_1 = 0. 
\end{split}
\end{equation}
Moreover, using $\varphi(v ) = v$ gives, following the computations in the previous sections, the equations
\begin{equation}
\label{eq:m1}
\begin{aligned}
&0&= \dfrac{\gamma \rho_1^2}{2} \left\{ \dfrac{\nu_1 + (1-p)\gamma^2}{\nu_1 + \gamma^2}\left( P(\rho_1) - (P(\rho_1)+(1-P(\rho_1))^2)m_1 \right)  + \dfrac{\gamma p}{\nu_1 + \gamma^2 }  (\bar v_1 - m_1) \right\}  \\
&&-\beta_1 (1-\rho_2)^{\alpha}\rho_1 m_1 + \beta_2 (1-\rho_1)^\alpha \rho_2m_2 
\end{aligned}
\end{equation}
\begin{equation}
\label{eq:m2}
\begin{aligned}
&0 &= \dfrac{\gamma \rho_2^2}{2} \left\{ \dfrac{\nu_2 + (1-p)\gamma^2}{\nu_2 + \gamma^2}\left( P(\rho_2) - (P(\rho_2)+(1-P(\rho_2))^2)m_2 \right)  + \dfrac{\gamma p}{\nu_2 + \gamma^2 }  (\bar v_2 - m_2) \right\}  \\
&&-\beta_2 (1-\rho_1)^{\alpha}\rho_2 m_2 + \beta_1 (1-\rho_2)^\alpha \rho_1m_1
\end{aligned}
\end{equation}
Alltogether we have obtained  3 equations for the unknowns $m_1,m_2,\rho_1, \rho_1$, that means we obtain a one-parameter family of solutions.
As the free parameter the collision invariant $\rho_1 +\rho_2=  \rho$, i.e. the total density,  is used. 
We proceed as follows:
defining  the monotone function $f(x) = \frac{(1-x)^\alpha}{x}, x \in [0,1]$ equation (\ref{eq:rho}) gives 

\[
\beta_2 f (\rho_1) = \beta_1 f(\rho_2)
\]
and we obtain 
the following equation for $\rho_2$
 \begin{equation}
 \label{eq:f}
 f^{-1} \left(\frac{\beta_1}{\beta_2}f(\rho_2)\right)+\rho_2 =  \rho.
 \end{equation}
 We observe that the function
 \[
 g(x) =  f^{-1} \left(\frac{\beta_1}{\beta_2}f(x)\right)+x
 \]
 is for any $\alpha >0$ monotone increasing for $x \in [0,1]$ with values in $[0,2]$.
 This  yields a unique solution $\rho_2 ( \rho)$ of equation (\ref{eq:f}) and then $\rho_1=\rho_1( \rho) = \rho - \rho_2( \rho)$.

Using these computations in equations (\ref{eq:m1}) and (\ref{eq:m2}) 
  yields, after solving a linear system, expressions for $m_1(\rho)$ and $m_2( \rho)$.
Using these expressions for $m_1(\rho)$ and $m_2( \rho)$
in the second equations of equations (\ref{eq:1}) and (\ref{eq:2})  and adding the two equations yields finally a dynamic equation for the sum of the densities

\begin{equation}
\label{eq:sum}
\partial_t  \rho + \partial_x ((m_1+m_2)( \rho)) = 0 .
\end{equation}

The densities  $\rho_1$ and $\rho_2$ on each of the lanes are then  given by equation
(\ref{eq:rho}) and $\rho_1 +\rho_2=  \rho$. .

\begin{remark}
In case $\alpha=1$ we obtain explicitly 
$f^{-1}(x) = \frac{1}{1+x}$ and
\[
\frac{1}{1+\frac{\beta_1}{\beta_2} (\frac{1}{\rho_2}-1)} +\rho_2 = \rho.
\]
Rewriting with $\beta= \frac{\beta_1}{\beta_2} $ we obtain the quadratic equation
\[\begin{aligned}
(1-\beta)\rho_2^2 + \rho_2 (1+\beta - \rho (1-\beta)) - \beta \rho =0
\end{aligned}
\]
with $\beta \in [0,\infty)$, $\bar \rho \in [0,2]$. We are looking for solutions  $\rho_2$ in $[0,1]$.
For $\beta =1 $ the solution is $\rho_2 = \frac{ \rho}{2}$.
For general $\beta$ one easily observes that there is a unique solution in $[0,1]$.
\end{remark}

%

\subsubsection{Slow switching regime}

In this regime we have $\beta_1,\beta_2 = O(\epsilon)$. Therefore two time scales are detectable, the one of interactions, fast in time, and the one of transport and reaction. 

In the limit $\epsilon\rightarrow 0^+$ we have that $\beta_1/\epsilon \rightarrow c_1$ and $\beta_2/\epsilon \rightarrow c_2$, where $c_1,c_2>0$. The fast and slow dynamics in equations \eqref{eq:multi_s1}-\eqref{eq:multi_s2} are
\[
\begin{cases}
&\partial_t \displaystyle \int_0^1 \varphi(v)f_1(x,v,t)dv = \dfrac{1}{\epsilon}\displaystyle \int_0^1 Q(f_1,f_1)(x,v,t)dv \\
&\partial_t \displaystyle\int_0^1 \varphi(v)f_1(x,v,t)dv + \partial_x\displaystyle \int_0^1 v\varphi(v)f_1(x,v,t)dv  \\
&\qquad\qquad =  -c_1(1-\rho_2 )^\alpha \displaystyle\int_0^1 \varphi(v)f_1(x,v,t)dv + c_2 (1-\rho_1)^\alpha \displaystyle \int_0^1 \varphi(v)f_2(x,v,t)dv
\end{cases}
\]
and
\[
\begin{cases}
&\partial_t \displaystyle \int_0^1 \varphi(v)f_2(x,v,t)dv = \dfrac{1}{\epsilon}\displaystyle \int_0^1 Q(f_2,f_2)(x,v,t)dv \\
&\partial_t \displaystyle\int_0^1 \varphi(v)f_2(x,v,t)dv + \partial_x\displaystyle \int_0^1 v\varphi(v)f_2(x,v,t)dv  \\
&\qquad\qquad =  -c_2(1-\rho_1 )^\alpha \displaystyle\int_0^1 \varphi(v)f_2(x,v,t)dv + c_1 (1-\rho_2)^\alpha \displaystyle \int_0^1 \varphi(v)f_1(x,v,t)dv. 
\end{cases}
\]

Once again, if $\epsilon \rightarrow 0^+$, we reach fast in time the local equilibrium and, in the quasi-invariant regime \eqref{eq:quasi-invariant} we may substitute the Boltzmann collision operator with one of Fokker-Planck-type, see Section \ref{sec:hom}, i.e.
\[
f_i(x,v,t) \approx \rho_i(x,t)f_i^\infty(v), \qquad i = 1,2. 
\]
We plug the local equilibria in the transport-reaction step to obtain a system of balance laws
\[
\begin{cases}
&\partial_t \left(\rho_1 \displaystyle \int_0^1 \varphi(v) f_1^\infty(v)dv \right) + \partial_x \left(\rho_1 \displaystyle \int_0^1 v\varphi(v)f_1^\infty(v) dv \right)  \\
&\qquad\qquad\qquad= -c_1 (1-\rho_2)^\alpha \rho_1 \displaystyle\int_0^1 \varphi(v)f_1^\infty(v)dv + c_2(1-\rho_1 )^\alpha \rho_2 \displaystyle\int_0^1 \varphi(v)f_2^\infty(v)dv.   \\
&\partial_t \left(\rho_2 \displaystyle \int_0^1 \varphi(v) f_2^\infty(v)dv \right) + \partial_x \left(\rho_2 \displaystyle \int_0^1 v\varphi(v)f_2^\infty dv \right)  \\
&\qquad\qquad\qquad= -c_2 (1-\rho_1)^\alpha \rho_2 \displaystyle\int_0^1 \varphi(v)f_2^\infty(v)dv + c_1(1-\rho_2 )^\alpha \rho_1 \displaystyle\int_0^1 \varphi(v)f_1^\infty(v)dv. 
\end{cases}\]

Since the only conserved quantity during interactions is the mass, we may set $\varphi(v) = 1$ and we obtain

\begin{equation}
\label{eq:hydro_slow}
\begin{cases}
\partial_t \rho_1 + \partial_x \left( \rho_1 \dfrac{P(\rho_1) + \frac{p}{\kappa_1} \bar v_1}{P(\rho_1) + (1-P(\rho_1))^2 + \frac{p}{\kappa_1}}\right) = -c_1 (1-\rho_2)^\alpha \rho_1  + c_2 (1-\rho_1)^\alpha \rho_2 \\
\partial_t \rho_2 + \partial_x \left( \rho_2 \dfrac{P(\rho_2) + \frac{p}{\kappa_2} \bar v_2}{P(\rho_2) + (1-P(\rho_2))^2 + \frac{p}{\kappa_2}} \right) = -c_2(1-\rho_1)^\alpha \rho_2  + c_1 (1-\rho_2)^\alpha \rho_1.
\end{cases}
\end{equation}

It is worth to remark that the introduced control modifies the definition of the macroscopic flux, as for the collision dominated case, whereas the switching dynamics is observable at the hydrodynamic level only in the case of slow switching. Lane changing, in particular, determines a reaction term in the system of first order macroscopic equations which disappears when looking at the evolution of $\rho_1 + \rho_2$.

\section{Numerical tests}\label{sec:numerics}

In this section we present several numerical tests which highlight the features of the proposed control strategy for traffic  dynamics in the multilane framework. In the collision dominated regime we obtained in Section \ref{sect:collision_dom} a decoupled set of first order macroscopic equations which is coherent with the tests performed in \cite{TZ1} for the one lane case. Therefore we will concentrate on the slow switching case.

First we will compare the introduced Boltzmann-type model with driver-assist controls with its derived hydrodynamic approximation.  In the following, we adopt a Direct Simulation Monte Carlo (DSMC) approach for the numerical solution of Boltzmann-type model, see \cite{PR,PT}. The evolution of the obtained macroscopic models are  computed through finite volume WENO schemes for conservation laws with non-convex fluxes \cite{QS,Shu}. In all the subsequent tests we will consider the function $P(\rho_i)$ defined in \eqref{eq:defP} with $\mu = 2$ and a recommended lane-dependent speed of the form $\bar{v}_i(\rho_i) = (1-\rho_i)$, $i = 1,2$. For the designed DSMC solver we consider $N = 2 \cdot10^{6}$ particles at time $t = 0$. The density is reconstructed on the phase space $(x,v) \in X\times V$, where $X =  [-2,2] $ and $V =  [0,1]$ by means of $N_x = 21$ gridpoints for the $x$ variable and $N_v = 128$ gridpoints for the $v$ variable. 

\subsection{Test 1}
In this test, we investigate numerically,  whether the derived hydrodynamics equations resulting from the kinetic model with  constrained binary interactions in the slow switching regime \ref{sect:slow} are consistent with the considered Boltzmann-type equations.

We briefly describe the numerical strategy adopted to solve the system of kinetic equations \eqref{eq:1_weak}-\eqref{eq:2_weak}. We begin by rewriting the considered inhomogeneous Boltzmann-type model in strong form for the slow switching regime. This read
\begin{equation}
\begin{cases}\vspace{0.25cm}
\partial_t f_1 + v\partial_x f_1 =& \dfrac{1}{\epsilon} Q(f_1,f_1)  - \beta_1(1-\rho_2)^\alpha f_1 + \beta_2(1-\rho_1)^\alpha f_2 \\
\partial_t f_2 + v\partial_x f_2 =& \dfrac{1}{\epsilon} Q(f_2,f_2)  - \beta_2(1-\rho_1)^\alpha f_2 + \beta_1(1-\rho_2)^\alpha f_1, 
\end{cases}
\end{equation}
where $Q(\cdot,\cdot)$ is a Boltzmann-type interaction operator defined in (\ref{collision}). 

 To describe the adopted numerical strategy we need to introduce a time discretization $t^n = n\Delta t$, $\Delta t>0$ and $n \in \mathbb N$. We adopt the following approach 

\begin{enumerate}
\item[$0)$] At time $t =0$ we consider the total initial density $f(x,v,0) = f_1(x,v,0)+ f_2(x,v,0)$ which is approximated by means of $N$ labeled particles $(x_i,v_i,\ell)$, $i = 1,\dots,N$ and $\ell = 1,2$ is the label indicating the initial lane. In this way the number of particles $(x_i,v_i,1)$ is $ \int_X \rho_1(x,0)dx$ and the number of particles $(x_i,v_i,2)$ is $\int_X \rho_2(x,0)dx$. 
\item[$1)$] In a single time step and in all space positions $x \in X$ and for a given $\epsilon>0$ we consider the problems
\begin{equation}
\partial_t \tilde f_i(x,v,t) = \dfrac{1}{\epsilon} Q(\tilde f_i,\tilde f_i), \qquad t \in (t^n, t^{n+1/3}], 
\; \tilde f_i(x,v,t^n) = f_i(x,v,t^n), \qquad i = 1,2,
\end{equation}
which are  solved through Nanbu's algorithm for Maxwellian molecules. We point the interested reader to \cite{PR,PT} for a general introduction and to \cite{TZ1} for a direct application to kinetic models for traffic dynamics. 
\item[$2)$] The transport steps are solved subsequently by considering 
\[
\partial_t \hat f_i(x,v,t) + v \partial_x \hat f_i(x,v,t) = 0, \qquad t \in (t^{n+1/3},t^{n + 2/3}], \;
\hat f_i(x,v,t^{n+1/3}) = \tilde f_i(x,v,t^{n+1/3}). 
\]
\item[$3)$] The lane changing dynamics is solved by considering
\[
\begin{cases}
\partial_t f_1(x,v,t) = -\beta_1 (1-\rho_2)^\alpha f_1 + \beta_2(1-\rho_1)^\alpha f_2, \qquad t \in (t^{n+2/3},t^{n + 1}]\\
\partial_t f_2(x,v,t) = -\beta_2 (1-\rho_1)^\alpha f_2 + \beta_1(1-\rho_2)^\alpha f_1,  \\
f_i(x,v,t^{n+2/3}) = \hat f_i(x,v,t^{n+2/3}). 
\end{cases}
\]
In particular, since the particles are labeled,  step (3) has a clear probabilistic interpretation. The lane $\ell = 1,2$ which counts for particles $(x_i,v_i,\ell)$ undergoes in $x_i \in X$ a loss of $\beta_\ell (1-\rho_m)$, $m \ne \ell$, particles and has a gain of $\beta_m(1-\rho_\ell)$ particles. 
\end{enumerate}

We consider the following initial condition
\begin{equation}
\begin{split}
f_1(x,v,0) &= \begin{cases}
2 & (x,v) \in [-1,0] \times [0,1/2]\\
0 & \textrm{otherwise},
\end{cases} \\
f_2(x,v,0) &= \begin{cases}
2 & (x,v) \in [-2,-1) \times [0,1/2]\\
0 & \textrm{otherwise},
\end{cases}
\end{split}
\label{eq:init_f12}
\end{equation}
mimicking the condition where vehicles have uniform velocities between $[0,1/2]$ and are uniformly distributed in $[-1,0)$ in the first lane, whereas vehicles of the second lane are uniformly distributed in $[-2,-1)$, i.e. the vehicles in the first lane precede the ones in the second lane.

 The initial conditions \eqref{eq:init_f12} translate to the following condition on the densities 
\begin{equation}
\label{eq:init_rho12}
\rho_1(x,0) = \chi(x \in [-1,0]), \qquad \rho_2(x,0) = \chi(x \in [-2,-1)), 
\end{equation}
being $\chi(\cdot)$ the indicator function
\[
\chi(x\in A) = \begin{cases}
1 & x \in A \\
0 & x \not\in A.
\end{cases}
\]

\begin{figure}
\centering
\subfigure[Lane 1]{
\includegraphics[scale = 0.35]{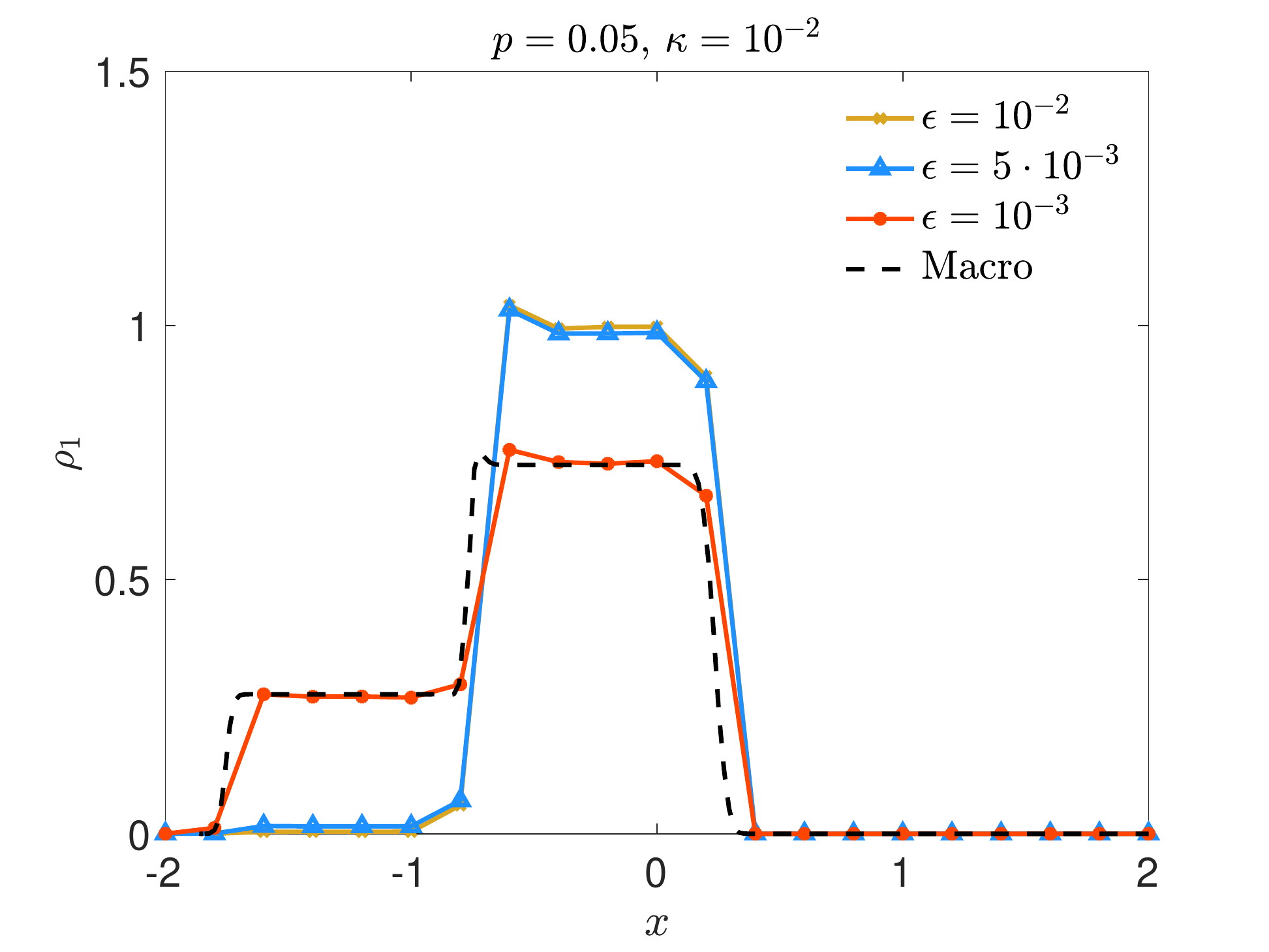}}
\subfigure[Lane 2]{
\includegraphics[scale = 0.35]{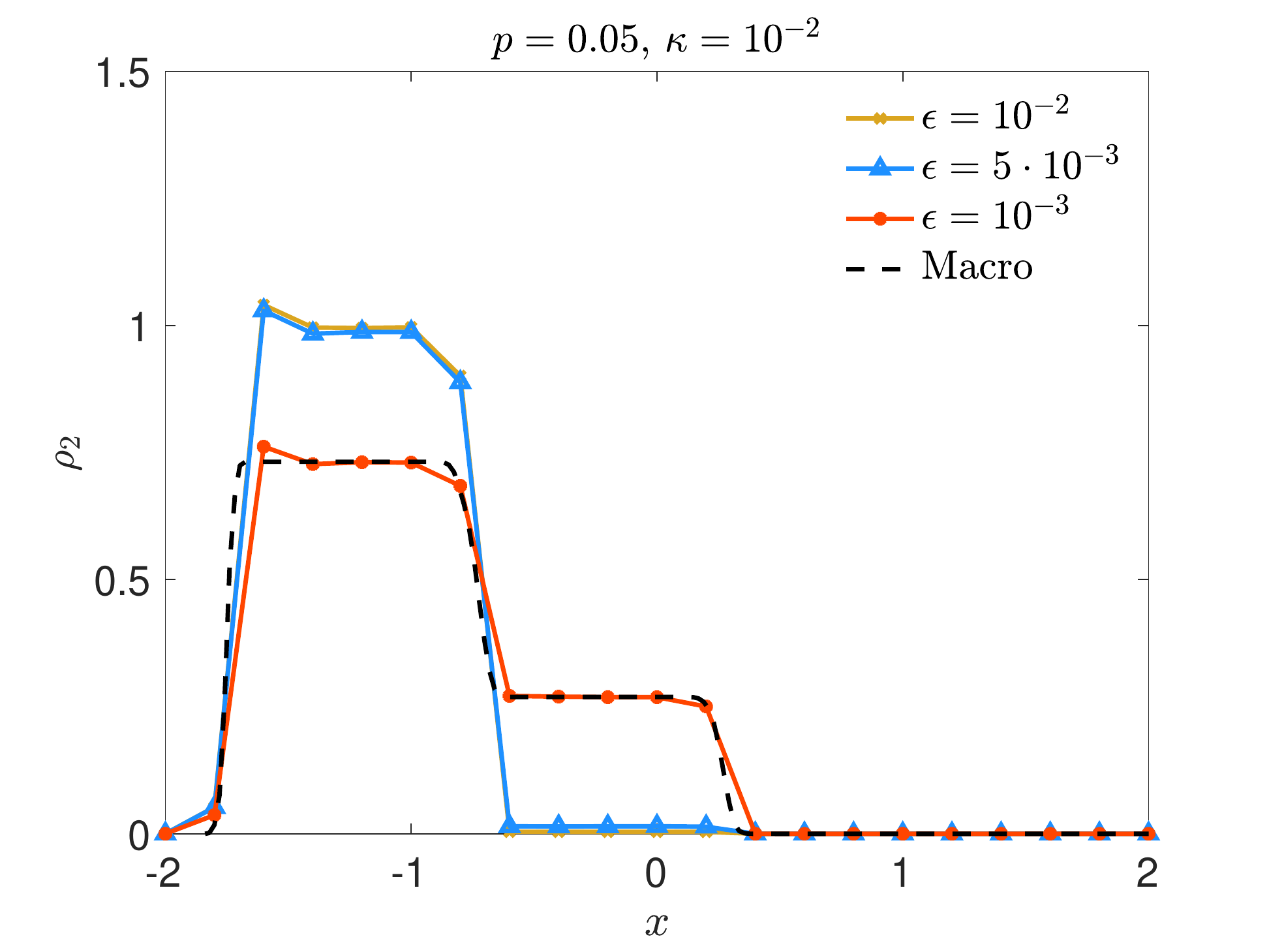}}
\caption{\textbf{Test 1}.  Evolution of the kinetic model \eqref{eq:1_weak}-\eqref{eq:2_weak} and of the derived hydrodynamic model \eqref{eq:hydro_slow} in the flow switching regime up to time $T = 0.2$ with penetration rate $p =0.05$ and penalization $\kappa = 10^{-2}$  with $\Delta t = \epsilon$ and decreasing $\epsilon = 10^{-2}, 5 \cdot 10^{-3}, 10^{-3}$ with an initial number of particles $N = 10^6$ in each lane. Other parameters are $\alpha = 2$, $\mu = 2$, $\beta_1 = \beta_2 = 2\epsilon$.  }
\label{fig:compare}
\end{figure}

In Figure \ref{fig:compare} we compare the densities $\rho_1,\rho_2$ computed from the numerical solution of the kinetic model \eqref{eq:1_weak}-\eqref{eq:2_weak} with the ones of the system of macroscopic equations \eqref{eq:hydro_slow} for several $\epsilon>0$. We can observe how the densities of the kinetic model converge towards the macroscopic ones for a sufficiently small $\epsilon$. 

\begin{figure}
\centering
\subfigure[Lane 1, $\kappa = +\infty$]{
\includegraphics[scale = 0.35]{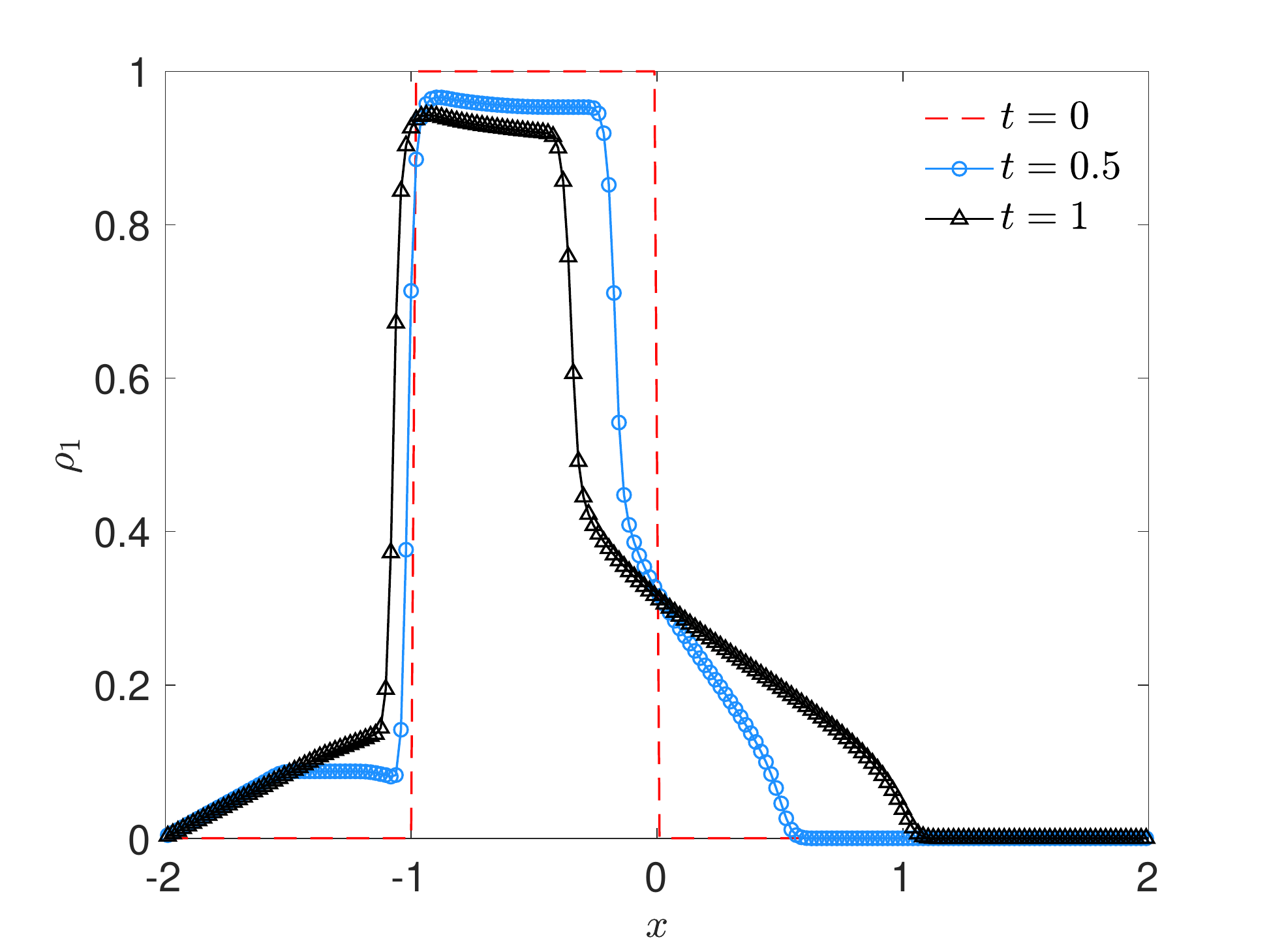}}
\subfigure[Lane 2, $\kappa = +\infty$]{
\includegraphics[scale = 0.35]{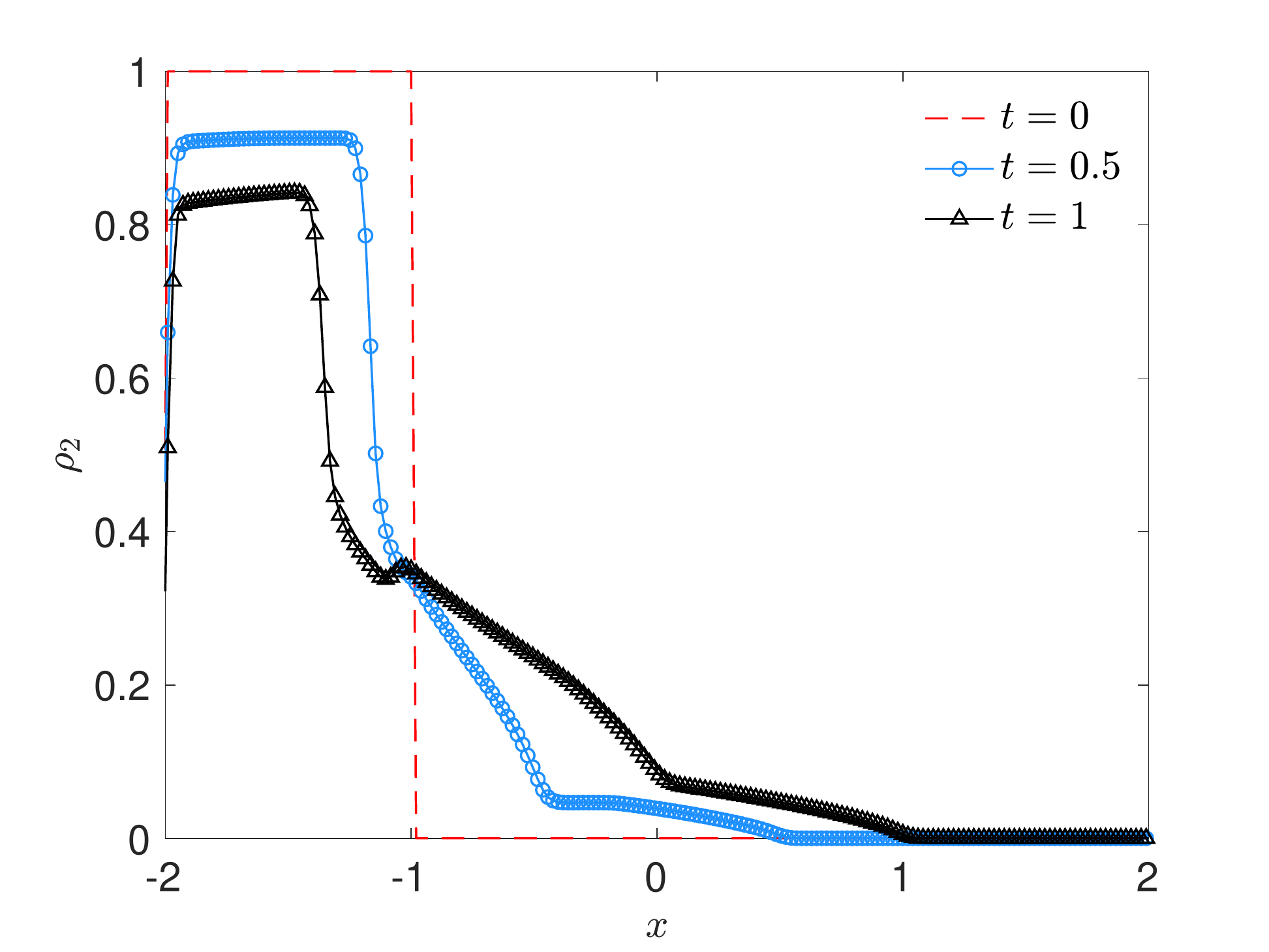} }\\
\subfigure[Lane 1, $p = 0.05,\kappa = 10^{-2}$]{
\includegraphics[scale = 0.35]{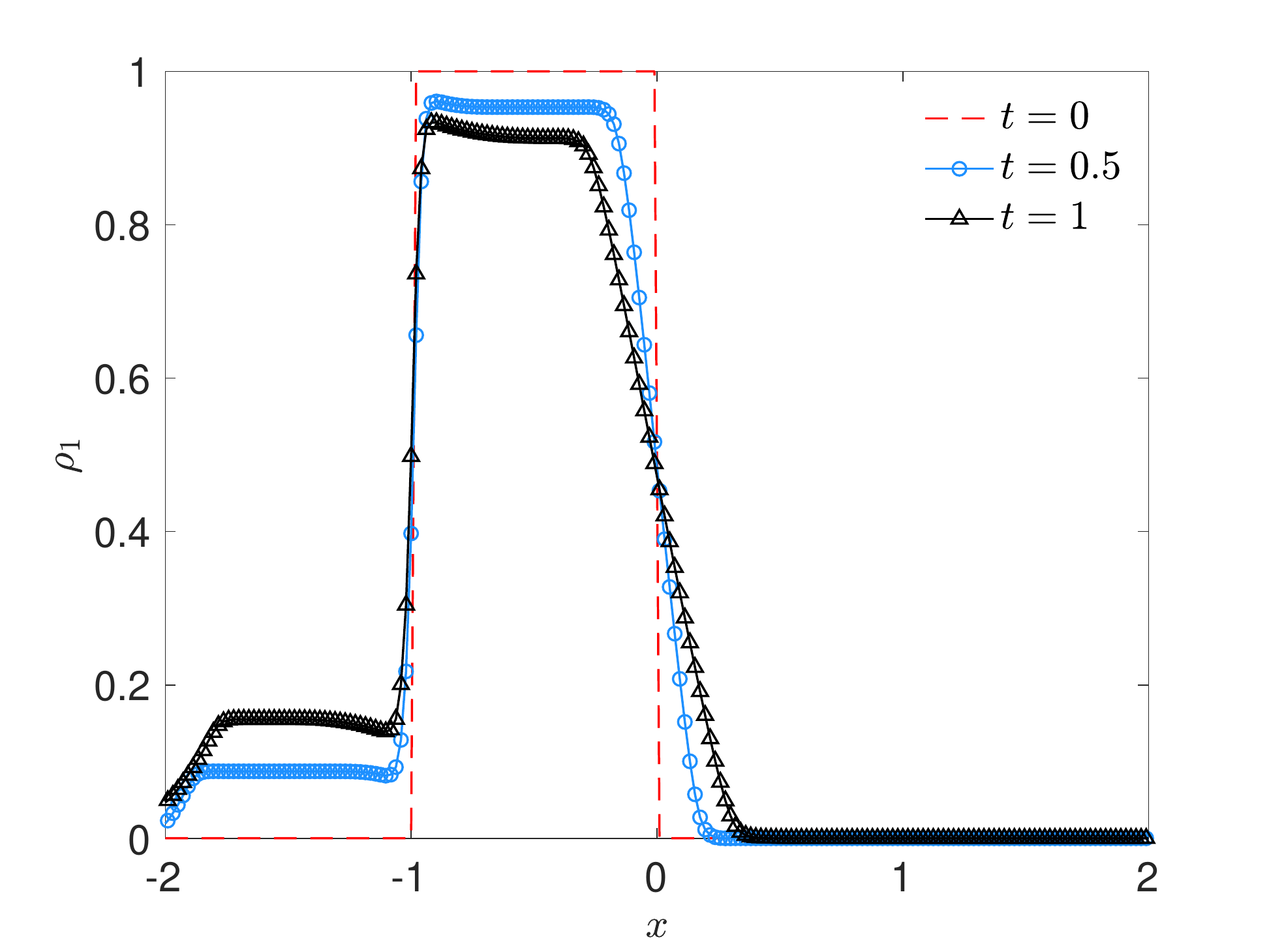}}
\subfigure[Lane 2, $p = 0.05, \kappa = 10^{-2}$]{
\includegraphics[scale = 0.35]{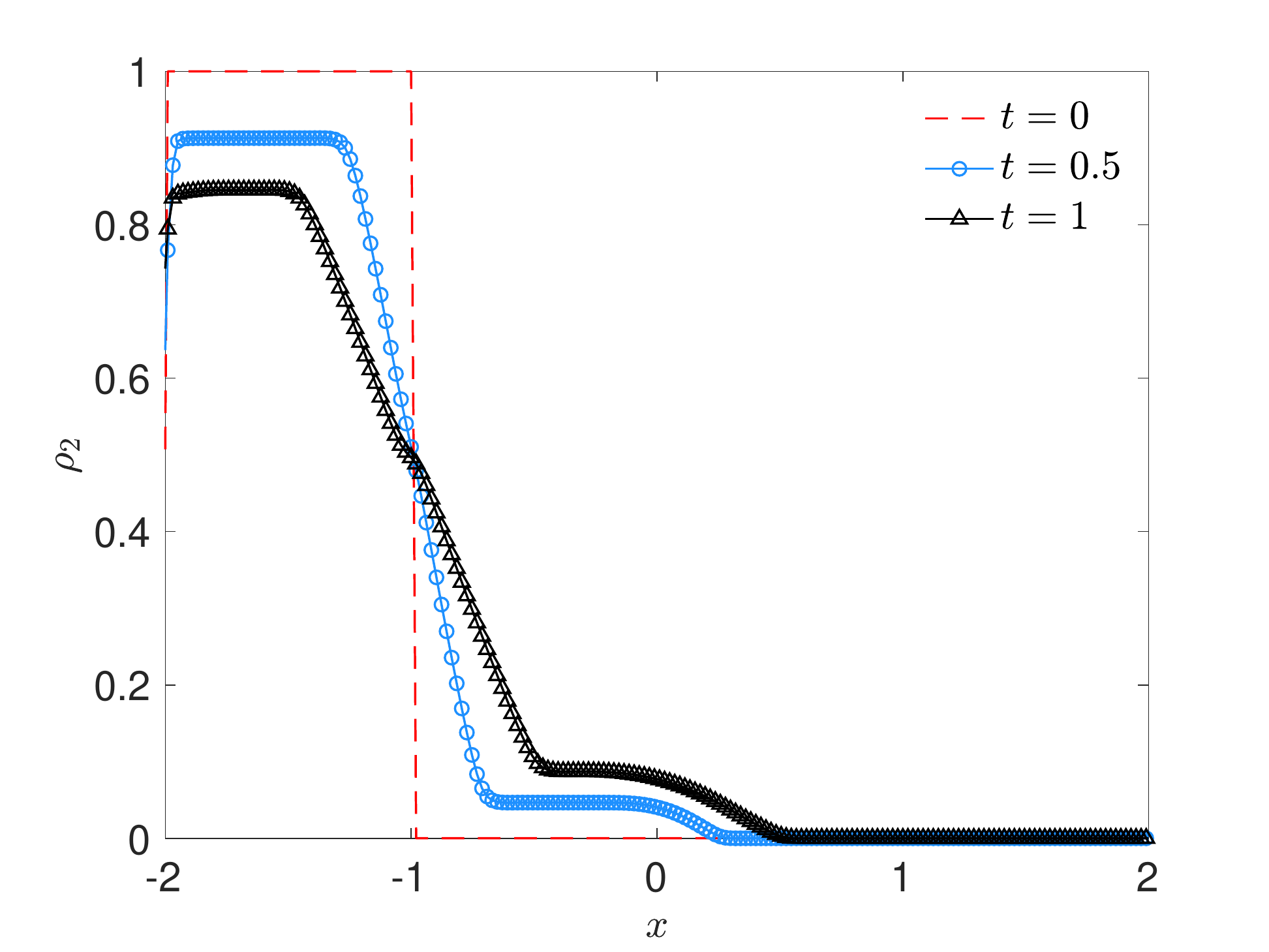}}
\caption{\textbf{Test 1}. Evolution at times $t = 0,1/2,1$ of the macroscopic model \eqref{eq:hydro_slow} in the uncontrained case $\kappa = +\infty$ (top row) and in the case with penetration rate $p = 0.05$ and $\kappa = 10^{-2}$. We considered rescaled switching parameters $c_1 = 0.1$, $c_2 = 0.2$}
\label{fig:macro_control}
\end{figure}

The action of the control is highlighted in Figure \ref{fig:macro_control} where in the top row we represent the evolution of $\rho_1$, $\rho_2$ at three different times $t = 0,\frac{1}{2},1$ in the unconstrained case. Due to the non-convexity  of the flux of \eqref{eq:hydro_slow} the solution of the Riemann problem determined by the initial condition \eqref{eq:init_rho12} is a combination of shocks and rarefactions. The introduced control has the effect of regularizing the flux such that for $\kappa \rightarrow 0^+$ converges to the solution of 
\[
\begin{cases}
\partial_t \rho_1 + \partial_x \left( \rho_1 \bar v_1(\rho_1)\right) = -c_1 (1-\rho_2)^\alpha \rho_1 + c_2(1-\rho_1)^\alpha \rho_2 \\
\partial_t \rho_2 + \partial_x \left( \rho_2 \bar v_2(\rho_2)\right) = -c_2 (1-\rho_1)^\alpha \rho_2 + c_1(1-\rho_2)^\alpha \rho_1, 
\end{cases}
\]
where in this test we considered $\bar v_i(\rho_i) = 1-\rho_i$. 

\subsection{Test 2}

\begin{figure}
\centering
\includegraphics[scale = 0.35]{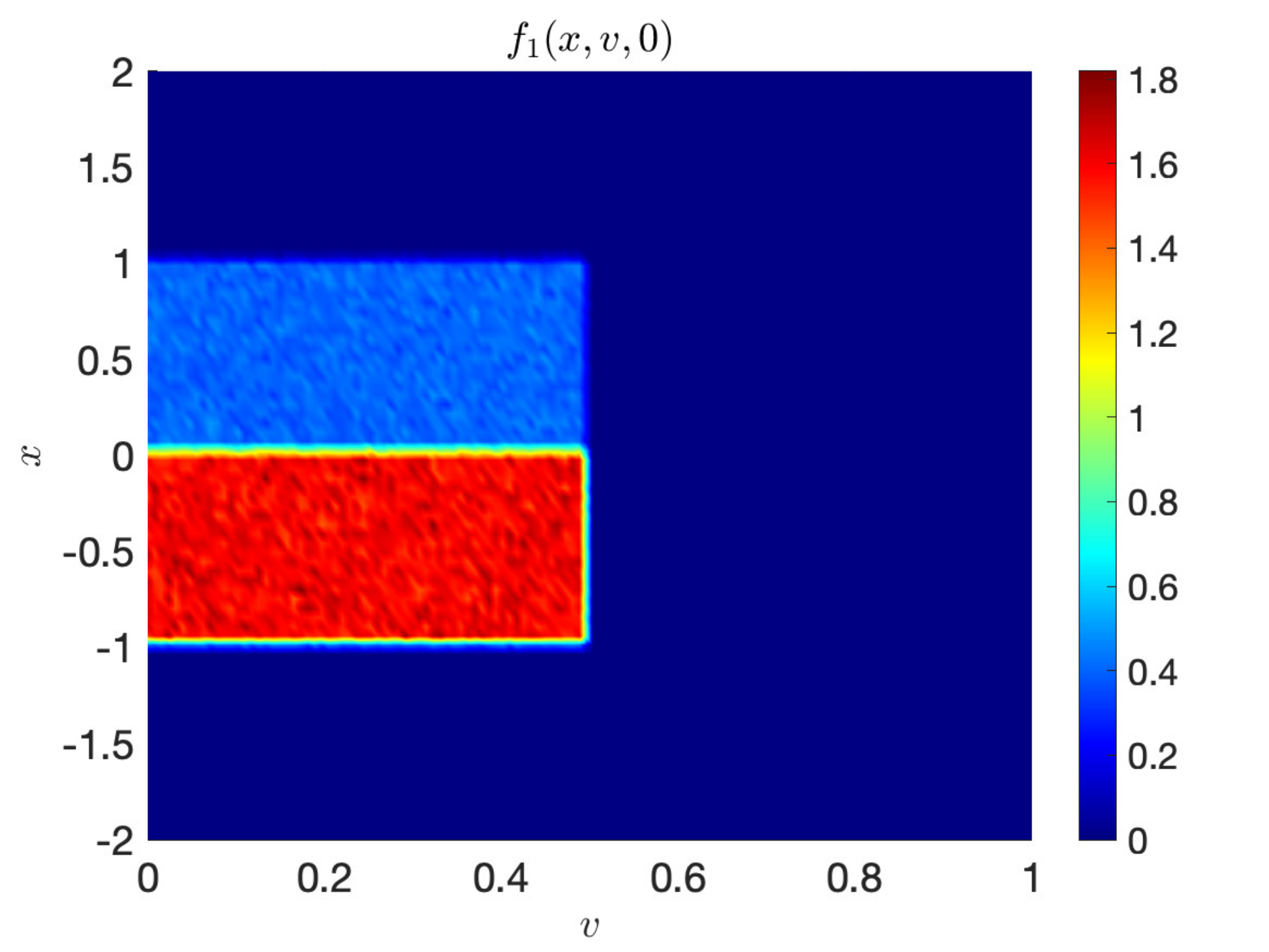}
\includegraphics[scale = 0.35]{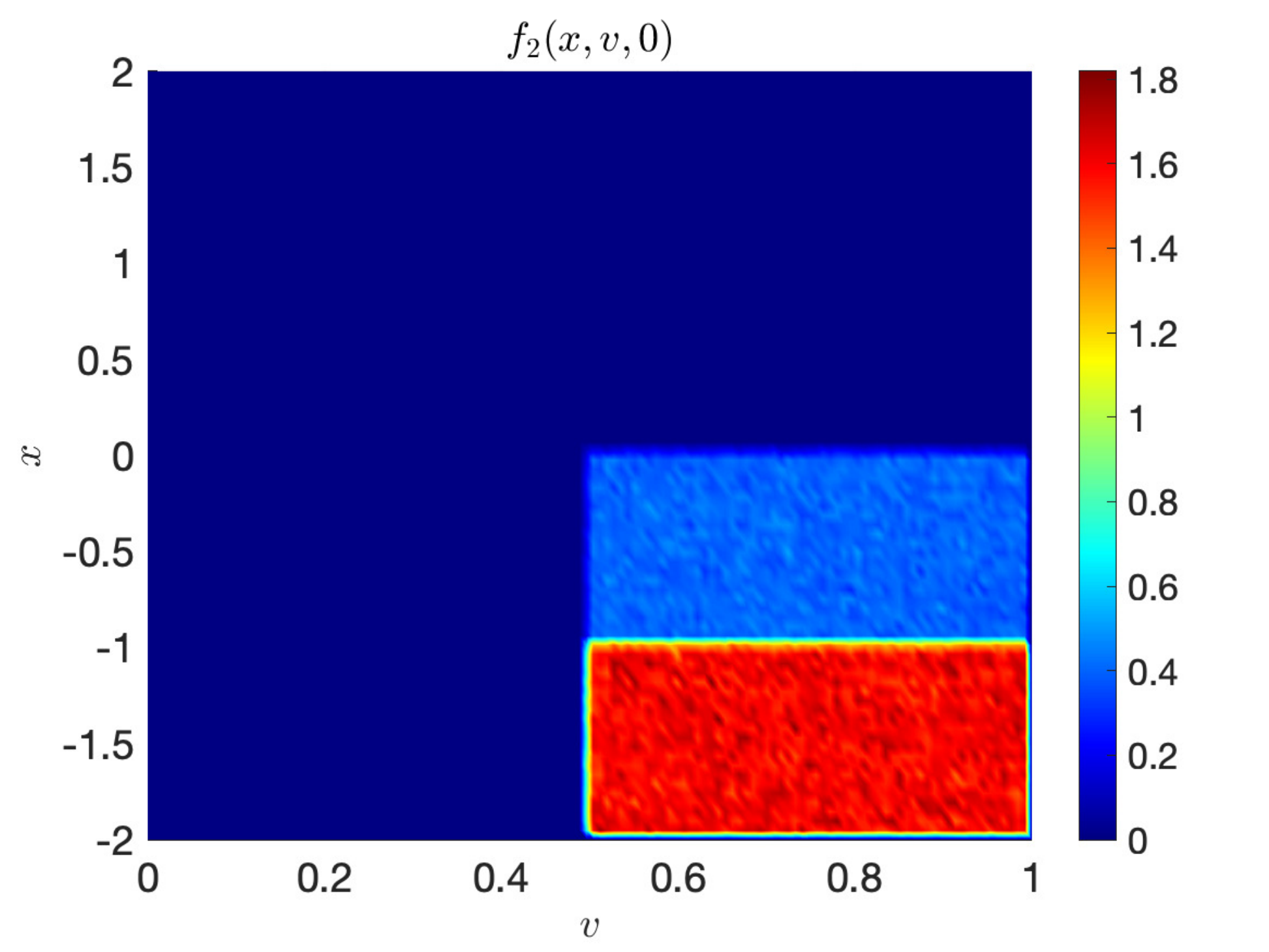} 
\caption{\textbf{Test 2}. Initial densities in the two lanes defined in \eqref{eq:init_f12_test2}. We considered $N = 10^6$ in each lane and the density is reconstructed with $N_v = 128$ gridpoints for the velocity variable and $N_x = 21$ gridpoints for the space veriable. }
\label{fig:init_test2}
\end{figure}

\begin{figure}
\centering
\includegraphics[scale = 0.35]{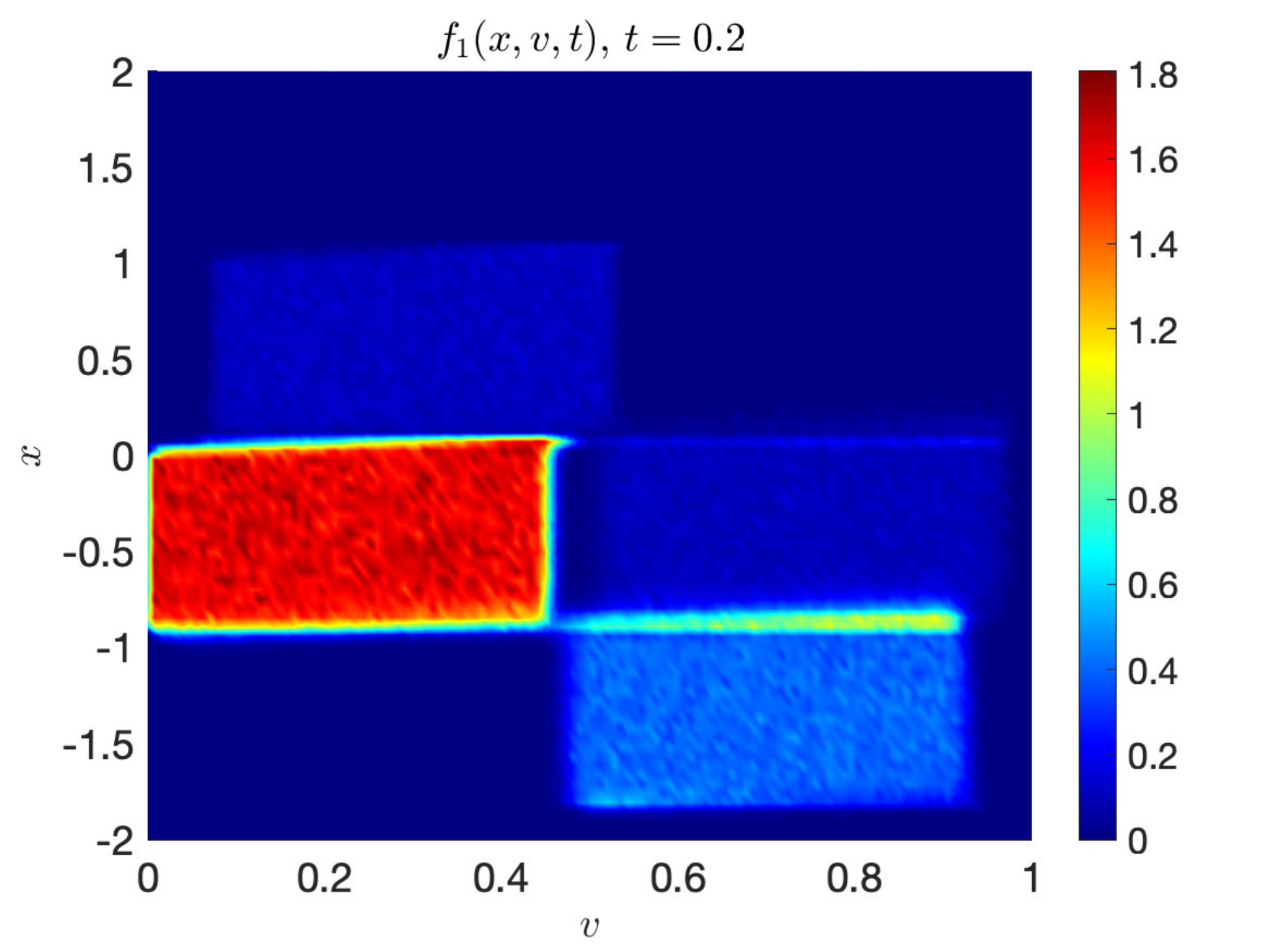}
\includegraphics[scale = 0.35]{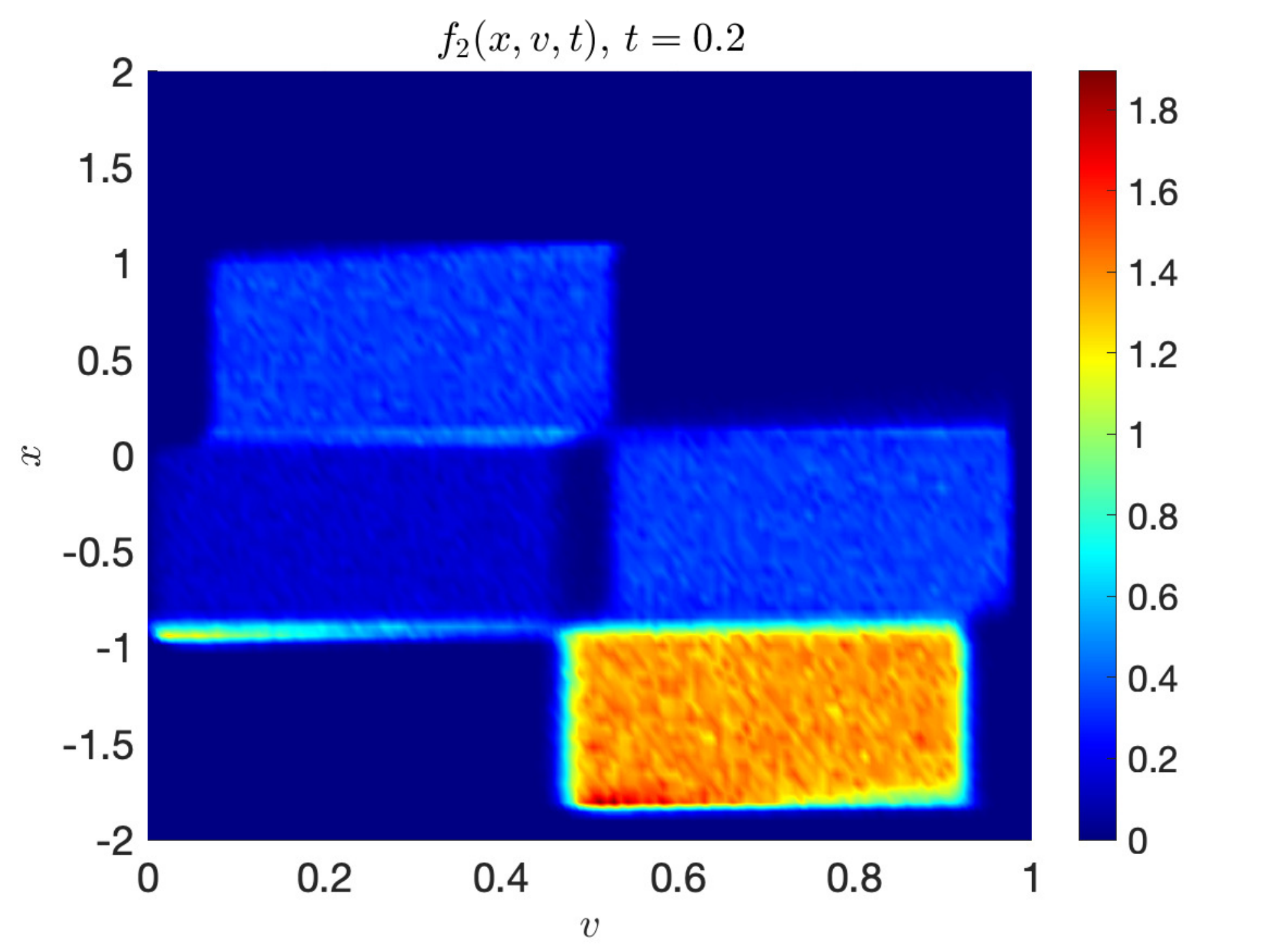} \\
\includegraphics[scale = 0.35]{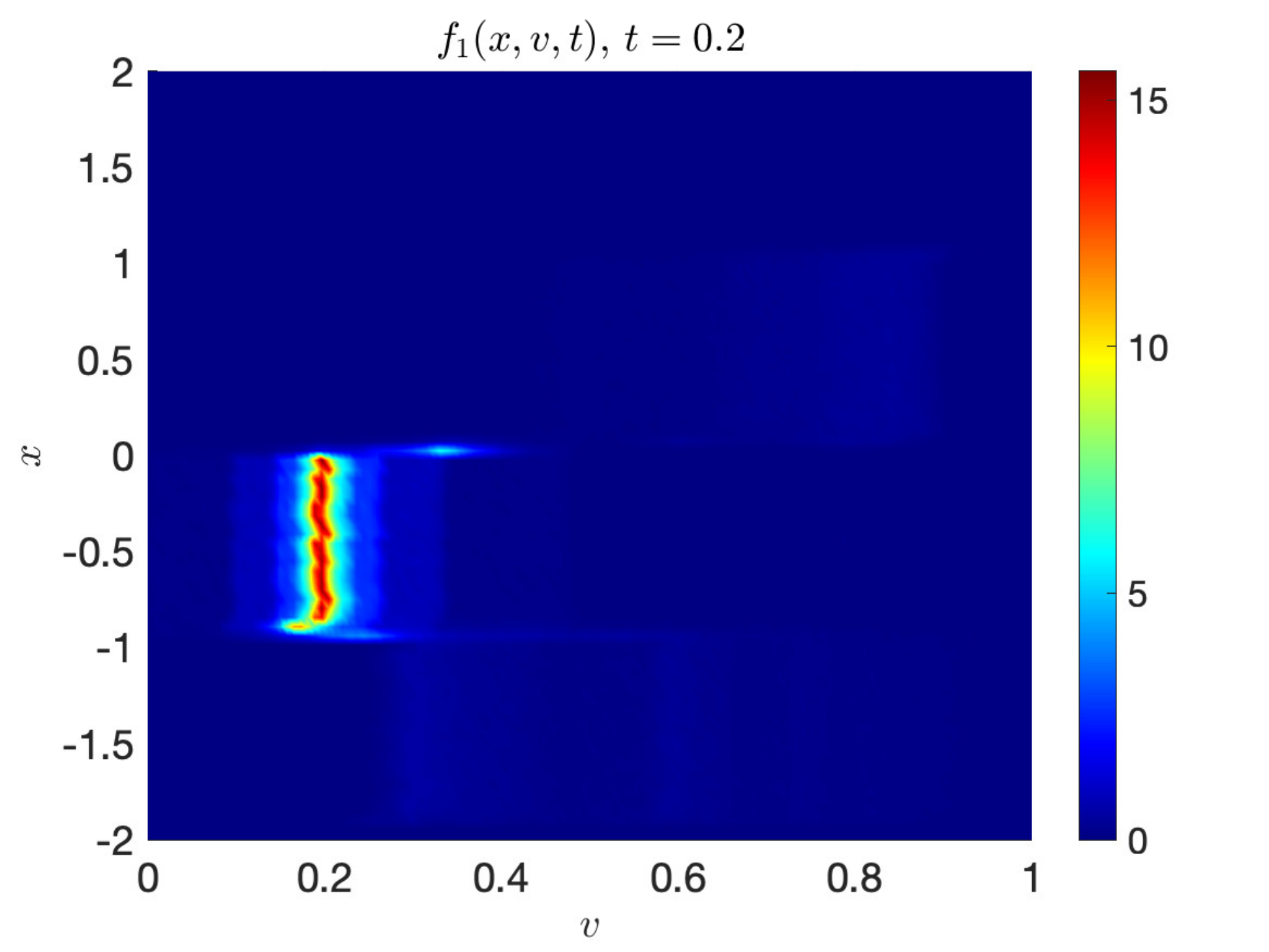}
\includegraphics[scale = 0.35]{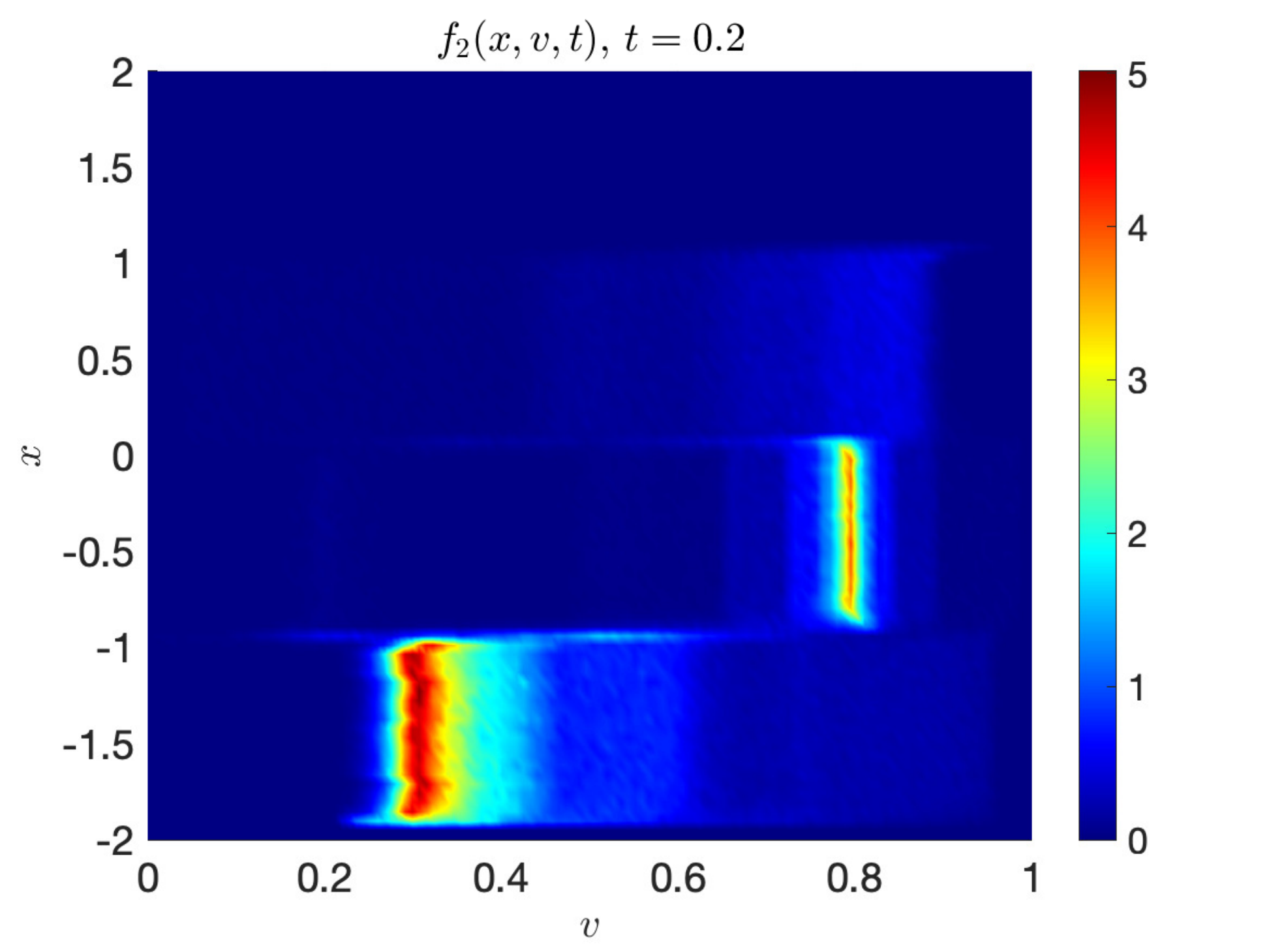}
\caption{\textbf{Test 2}. Densities of the first lane (left) and of the second lane (right) in the non controlled case (top row) and in the controlled case (bottom row). We considered $p = 0.05$ and $\kappa = 10^{-2}$. The switching coefficient $\beta_i$ are defined as in \eqref{eq:beta_m12} with $a = 0.2$.    }
\label{fig:test2_evo}
\end{figure}

In this test we exploit the structure of the introduced kinetic model to investigate the influence of the mean speed in lane changing dynamics. Indeed, the possibility of lane switching may depend also on the difference on mean speed between lanes such that radical difference of speed deters the drivers in changing lane.   We can take into account the described behavior by defining the following time dependent switching parameters
\begin{equation}
\label{eq:beta_m12}
\beta_i(x,t) = \dfrac{\epsilon}{|m_1(x,t) - m_2(x,t)| + a}, \qquad i = 1,2,
\end{equation}
where $a\in (0,1)$ and 
\[
\rho_i(x,t)m_i(x,t) = \int_0^1 vf_i(x,v,t)dv, \qquad i = 1,2. 
\]
In the following, we will consider the following initial condition 
\begin{equation}
\begin{split}
f_1(x,v,0) = \begin{cases}
2 \rho_L & (x,v) \in [-1,0) \times [0,1/2] \\
2 \rho_R & (x,v) \in [0,1] \times [0,1/2] \\
0 & \textrm{otherwise},
\end{cases} \\
f_2(x,v,0) = \begin{cases}
2 \rho_L & (x,v) \in [-2,-1) \times [1/2,1] \\
2 \rho_R & (x,v) \in [-1,0] \times [1/2,1] \\
0 & \textrm{otherwise},
\end{cases} 
\end{split}
\label{eq:init_f12_test2}
\end{equation}
with $\rho_L= 0.8$ and $\rho_R = 0.2$. The condition \eqref{eq:init_f12_test2} corresponds at the macroscopic level to consider 
\begin{equation}
\rho_1(x,0) = 
\begin{cases}
\rho_L & x \in [-1,0)\\
\rho_R & x \in [0,1] \\
0 & \textrm{otherwise}, 
\qquad
\end{cases}
\rho_2(x,0) = 
\begin{cases}
\rho_L & x \in [-2,-1)\\
\rho_R & x \in [-1,0] \\
0 & \textrm{otherwise}, 
\end{cases}
\end{equation}
see Figure \ref{fig:init_test2}. In Figure \ref{fig:test2_evo} we represent the evolution of the densities $f_1,f_2$ in the unconstrained case (top row) and in the case with driver-assist controls. As before we considered a penetration rate $p = 0.05$ and penalization term $\kappa = 10^{-2}$. The controls steer the velocities towards the desired speeds  $\bar v_i = 1-\rho_i$. 

The $\beta_i$ defined in \eqref{eq:beta_m12} induces stronger lane changing if vehicles' speed are similar, i.e. in our case in $x \in [-1,0]$, whereas inhibits switching for larger  differences of the speeds. This behavior can be observed in the first row of Figure \ref{fig:test2_evo} where the kinetic model is evolved up to time $t = 0.2$ using initial condition \eqref{eq:init_f12_test2}. In the second row we show the effect of the considered control in aligning the velocities of the vehicles. 

\begin{figure}
\includegraphics[scale = 0.4]{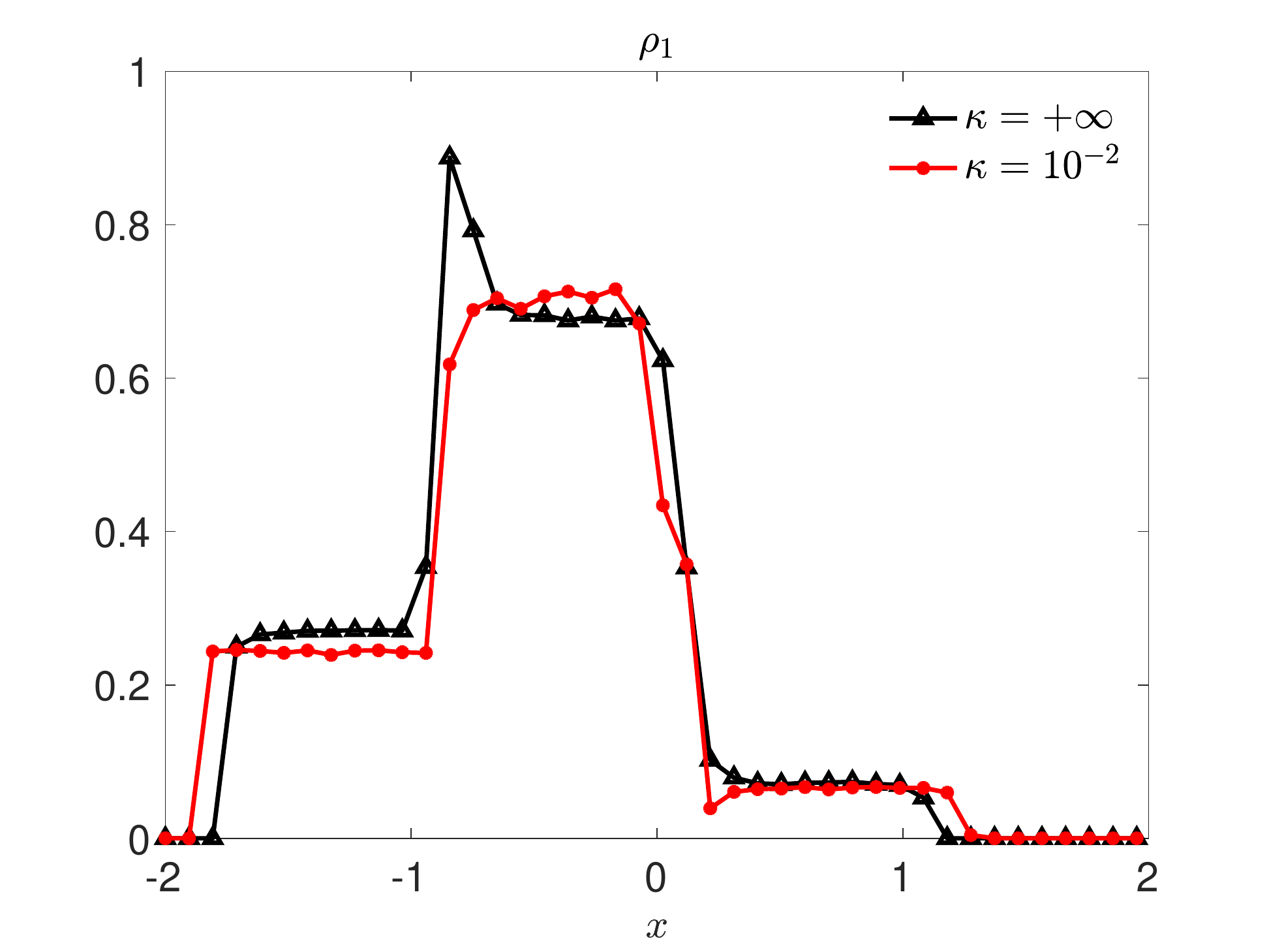}
\includegraphics[scale = 0.4]{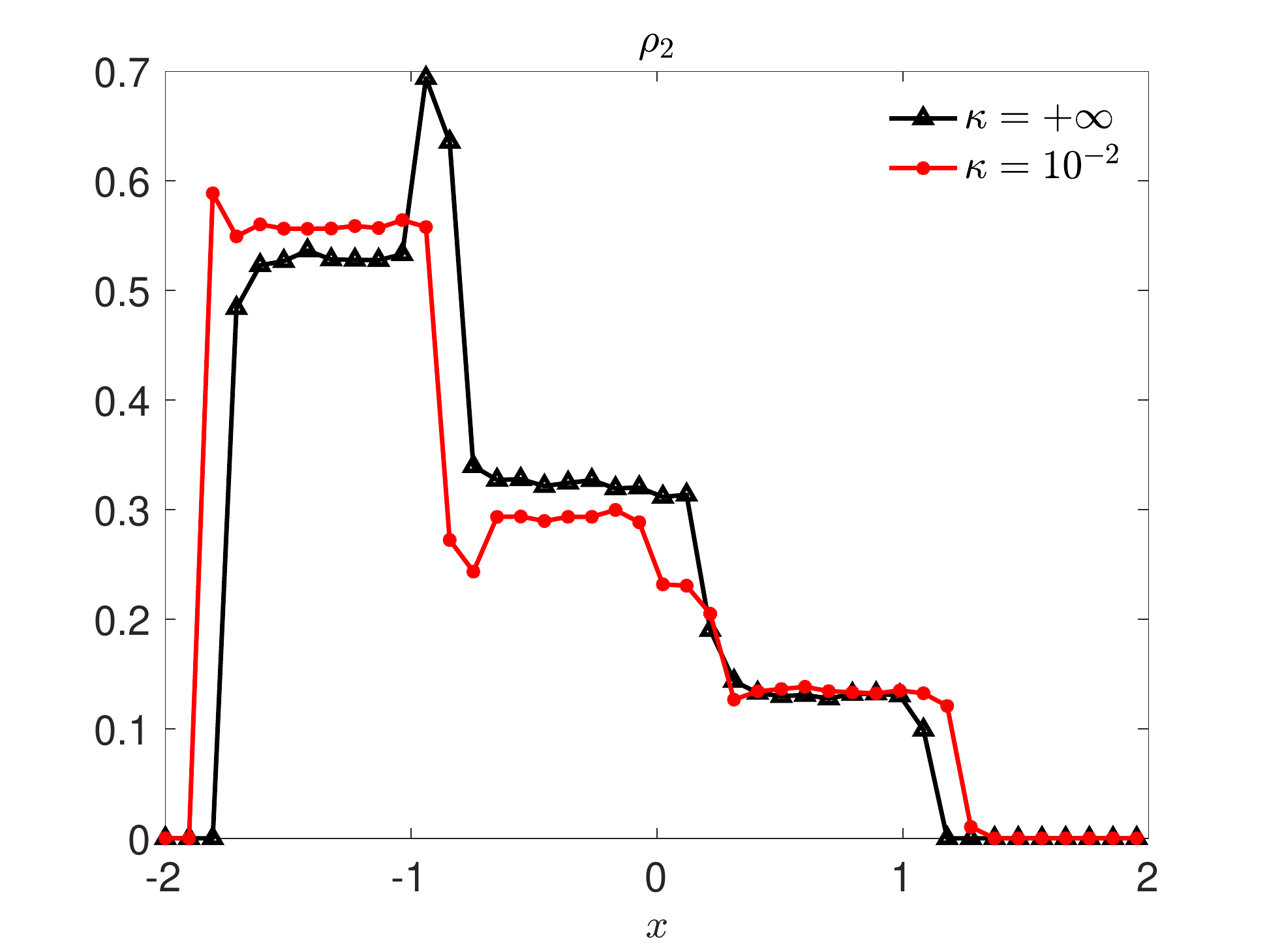}
\caption{\textbf{Test 2}. Macroscopic densities of the first and second lane in the controlled and non controlled settings. The densities are obtained through integration of the numerical solutions of the kinetic model with switching parameters in \eqref{eq:beta_m12}. }
\label{fig:test2_macro}
\end{figure}

In Figure \ref{fig:test2_macro} we represent $\rho_1$ and $\rho_2$ at time $t = 0.2$ obtained by direct integration of the kinetic densities $f_1(x,v,t)$ and $f_2(x,v,t)$ in Figure \ref{fig:test2_evo} in both the non controlled case, i.e. $\kappa = +\infty$, and in the controlled framework with $\kappa = 10^{-2}$ for a global penetration rate $p = 0.05$. The action of the control is clearly observable through a regularization of the flux near discontinuities. 

\section*{Conclusion}
In this paper we developed a hierarchical approach to the control of multilane traffic dynamics to understand the global impact of automated vehicles on the flow. In agreement with recent results on kinetic modelling for automated traffic, we consider anisotropic binary interactions between vehicles which include probabilistically the presence of driver-assist vehicles in the traffic flow. In addition,  lane changing operators are considered at the kinetic level and we showed the robustness of the controls with respect to the considered switching mechanism in terms of the alignment of velocities. Finally, fluid models are derived from the space inhomogeneous setting. To this end we considered a suitable closure method which rely on the possibility of the introduced kinetic model to provide explicit information on the constrained  speed distribution at equilibrium. The resulting fluid models are original macroscopic traffic models synthesizing controls embedded at the microscopic scale and that depend on the lane switching frequency.

\section*{Acknowledgements}
M. Z. is member of  GNFM (Gruppo Nazionale per la Fisica Matematica) of INdAM, Italy.
M. Z. acknowledges fellowship provided by Hausdorff Institute for Mathematics, Bonn, Germany, for the Junior Trimester on Kinetic Theory.

\end{document}